\documentclass[12pt,a4paper]{article}

\usepackage[usenames]{color}
\usepackage{amsmath,amsfonts}
\usepackage{mathrsfs}

\newtheorem{lemma}{Lemma}
\newtheorem{proposition}{Proposition}

\let\scr\mathscr

\def \Limsup{\mathop{\overline{\lim}}\limits}

\def\Pb{\mathbf{P}}

\def\Ex{\mathbf{E}}

\def\EE{\mathbb{E}}

\def\RR{\mathbb{R}}

\def\sgn{{\rm sgn}}
\def\1{\mbox{1\hspace{-.25em}I}}

\begin{document}
\title{On  Score-Functions   and    Goodness-of-Fit Tests  for
  Stochastic Processes}
\author{Yu.A. Kutoyants\\
{\small Laboratoire de Statistique et Processus, Universit\'e du Maine}\\
 Le Mans,  France\\
{\small and}\\
{\small  Laboratory of Quantitive Finance, Higher School of Economics}\\
{  Moscow, Russia}\\
}

\date{}

\maketitle
\begin{abstract}
The problems of the construction of the asymptotically distribution free
goodness-of-fit
 tests for three models of stochastic processes  are considered. The null
 hypothesis for all models is composite parametric.  All tests are
 based on the  score-function processes, where the unknown parameter is replaced
 by the MLE. We show that a special  change of time  transforms the
 limit score-function processes into  the Brownian bridge. This property allows us to construct
 the asymptotically distribution free tests for the following three models of
 stochastic processes : dynamical systems with small noise, ergodic diffusion
 processes, inhomogeneous Poisson processes  and nonlinear AR time series.
\end{abstract}
\noindent MSC 2000 Classification: 62M02,  62G10, 62G20.

\bigskip
\noindent {\sl Key words}: \textsl{
 Cram\'er-von Mises type
  tests,  dynamical systems, small noise, ergodic diffusion process,
  inhomogeneous Poisson processes, nonlinear AR, goodness-of-fit tests, asymptotically
  distribution free tests.}

\section{Introduction}

We consider the problem of the construction of  asymptotically distribution
free goodness-of-fit tests for the three models of stochastic processes
observed in continuous time: small
noise diffusion, ergodic diffusion and inhomogeneous Poisson process.  We assume that
under the basic hypotheses the models depend  on some unknown one-dimensional
parameter.

Let us recall what happens in the similar problem in the well-known
i.i.d. model. Suppose that  we observe $n$ i.i.d. r.v.'s
$\left(X_1,\ldots,X_n\right)=X^n$
 with continuous distribution function
$F\left(x\right)$ and the basic (null)  hypothesis
  is parametric
\begin{align*}
{\cal H}_0\quad :\qquad \qquad F\left(x\right)=F\left(\vartheta
,x\right),\quad \vartheta \in \Theta
\end{align*}
where $F\left(\vartheta
,x\right)$ is known smooth function of $\vartheta\in \Theta =\left(a ,b \right) $ and $x$.

We have to construct a goodness-of-fit (GoF) test $\hat\psi_n$ which belongs to
the class  ${\cal K}_\alpha $ of tests of asymptotic size $\alpha $, i.e.,
\begin{align*}
{\cal K}_\alpha =\left\{\bar\psi _n\quad :\quad \Ex_\vartheta \bar\psi
_n=\alpha +o\left(1\right)\right\} \qquad {\rm for\;\; all}\quad \vartheta \in
\Theta .
\end{align*}

Introduce the Cram\'er-von Mises type statistic
\begin{align*}
\delta _n=n\int_{-\infty }^{\infty }\left[\hat
  F_n\left(x\right)-F\left(\hat\vartheta _n,x\right)\right]^2 {\rm
  d}F\left(\hat\vartheta _n,x\right) ,\quad \hat
  F_n\left(x\right)=\frac{1}{n}\sum_{j=1}^{n}\1_{\left\{X_j<x\right\}},
\end{align*}
where $\hat\vartheta _n$ is the maximum likelihood estimator (MLE) and $\hat
  F_n\left(x\right)$ is the empirical distribution function.

Note that if $\Theta =\left\{\vartheta _0\right\}$ (simple basic hypothesis),
then
\begin{align*}
\delta _n&=n\int_{-\infty }^{\infty }\left[\hat
  F_n\left(x\right)-F\left(\vartheta _0,x\right)\right]^2 {\rm
  d}F\left(\vartheta _0,x\right) \\
&\qquad \Longrightarrow
\int_{-\infty }^{\infty }B\left(F\left(\vartheta _0,x\right)\right)^2{\rm
  d}F\left(\vartheta _0,x\right)=\int_{0}^{1}B\left(s\right)^2{\rm d}s\equiv \Delta ,
\end{align*}
where $s=F_0\left(\vartheta ,x\right)$ and $B\left(s\right),0\leq s\leq 1$ is
a Brownian bridge. Therefore the test
$\hat\psi _n=\1_{\left\{\delta _n>c_\alpha \right\}}$ where $c_\alpha $ is the
solution of  equation $\Pb\left(\Delta >c_\alpha \right)=\alpha $ belongs
to ${\cal K}_\alpha $. Moreover  it is {\it asymptotically distribution free}
(ADF), because the limit distribution of the statistic $\delta _n$ does not
depend on $F\left(\vartheta_0 ,\cdot \right)$.

Let us return to the parametric basic hypothesis and  suppose that the
  model is sufficiently regular to satisfy the presented below expansion of
  the MLE:
\begin{align*}
&u_n\left(x\right)=\sqrt{n}\left(\hat F_n\left(x\right)-F\left(\hat\vartheta
_n,x\right)\right)\\
&\quad=\sqrt{n}\left(\hat F_n\left(x\right)-F\left(\vartheta
,x\right)\right)+\sqrt{n}\left( F\left(\vartheta,x\right)-F\left(\hat\vartheta
_n,x\right)\right)\\
&=B_n\left(x\right)-\sqrt{n}\left(\hat\vartheta _n-\vartheta
  \right)\,\dot F\left(\vartheta,x \right)+o\left(1\right).
\end{align*}
Here $\dot F\left(\vartheta,x \right)$ means the derivative of
$F\left(\vartheta ,x\right)$ w.r.t. $\vartheta $.  The first term
$B_n\left(x\right)=\sqrt{n}\left(\hat F_n\left(x\right)-F\left(\vartheta
,x\right)\right)$ as before converges to the Brownian bridge
$B\left(F\left(\vartheta ,x\right)\right) $ and the MLE admits the
representation
\begin{align*}
\sqrt{n}\left(\hat\vartheta _n-\vartheta \right)
=\frac{1}{\sqrt{n}}\sum_{j=1}^ {n}\frac{\dot \ell\left(\vartheta
  ,X_j\right)}{{\rm I}\left(\vartheta \right)} +o\left(1\right)=\int_{}^
{}\frac{\dot \ell\left(\vartheta ,y\right)}{{\rm I}\left(\vartheta
  \right)}{\rm d}B_n\left(y\right) +o\left(1\right).
\end{align*}
Here $\ell\left(\vartheta
  ,x\right)=\ln f\left(\vartheta
  ,x\right) $, $f\left(\vartheta
  ,x\right) $ is the density function  and   $ {\rm
    I}\left(\vartheta \right)$ is the Fisher information. It can be shown that
\begin{align}
u_n\left(x\right)&\Longrightarrow
B\left(F\left(\vartheta
,x\right)\right)-\int_{}^{ }\frac{\dot
  \ell\left(\vartheta ,y\right)}{\sqrt{{\rm I}\left(\vartheta
    \right)}}\;{\rm
  d}B\left(F\left(\vartheta,y\right)\right)\; \int_{-\infty }^{x }\frac{\dot
    \ell\left(\vartheta
  ,y\right)}{ \sqrt{{\rm I}\left(\vartheta     \right)}}\; {\rm
  d} F\left(\vartheta,y\right)\nonumber\\
&\quad=B\left(s\right)-\int_{0}^{1}h\left(\vartheta,v\right)\,{\rm
  d}B\left(v\right)\int_{0}^{s}h\left(\vartheta,v\right)\,{\rm d}v\equiv
u\left(s\right),
\label{0}
\end{align}
where $s=F\left(\vartheta ,x\right) $,
\begin{align*}
 h\left(\vartheta,s\right)=\frac{\dot \ell\left(\vartheta
  , F^{-1}_\vartheta \left(s\right)\right)}{ \sqrt{{\rm I}\left(\vartheta
    \right)}}, \qquad \int_{0}^{1}h\left(\vartheta,v\right)^2\,{\rm d}v=1.
\end{align*}
Therefore  $u_n\left(\cdot \right) $ converges to the random function $u\left(\cdot
\right)$
and this allows us  to prove (see Darling \cite{Dar}) the convergence
$$
\delta _n\Longrightarrow \int_{0}^{1}u\left(s\right)^2{\rm d}s.
$$
Hence the test based on $\delta _n$ is not ADF because the limit
distribution of the statistic $\delta _n$ depends on $F\left(\vartheta
,x\right)$. This makes the choice of the threshld $c_\alpha $ a more difficult
problem.

One possibility to obtain ADF test is to find a  linear transformation of
$u\left(\cdot \right)$ into Wiener process:
$
L_W\left[u\right]\left(s\right)=w\left(s\right).$
Then
$$
\int_{-\infty }^{\infty }\Bigl(L_W\left[u\right]\left(F\left(\vartheta
,x\right)\right)\Bigr)^2{\rm d} F\left(\vartheta
,x\right)=\int_{0}^{1}w\left(s\right)^2{\rm d}s \equiv \hat\delta .
$$
Therefore if we take the  statistics
\begin{align*}
\hat\delta _n=\int_{-\infty }^{\infty
}\Bigl(L_W\left[u_n\right]\left(F(\hat\vartheta _n,x) \right)\Bigr)^2{\rm d}F(\hat\vartheta _n,x)
\end{align*}
and verify the convergence $\hat\delta _n\Rightarrow \hat\delta $, then the test
$
\hat\psi _n=\1_{\left\{\hat\delta_n>d_\alpha  \right\}}$  with
 $\Pb\left(\hat\delta >d_\alpha \right)=\alpha
$
 is ADF  and belongs to ${\cal K}_\alpha $. Note that such transformation
$L_W\left[u \right]$ was proposed by Khmaladze
\cite{Kh81} (see also the different proof of it   in \cite{KK}).

In the present work we consider a similar problem of construction of ADF GoF tests
for stochastic processes, for which we suggest a much simpler transformation of the corresponding limit
statistics into the Brownian bridge.

The goal of this work is to study the GoF tests for three models of
observations of continuous time stochastic processes:  diffusion processes
$X^\varepsilon =\left(X_t,0\leq t\leq T\right)$
with small diffusion coefficient ($\varepsilon \rightarrow 0$),  ergodic
diffusion processes $X^T
=\left(X_t,0\leq t\leq T\right)$, $T\rightarrow \infty $ and  $\tau _*
$-periodic Poisson processes $X^n =\left(X_t,0\leq t\leq T=\tau _* n\right)$,
$n\rightarrow \infty $.  For all three models we introduce the corresponding
score-function processes (SFP) $U_\varepsilon\left(\cdot \right)
,U_T\left(\cdot \right)$ and $U_n\left(\cdot \right)$ and then we show that
the Cram\'er-von Mises type statistics based on these SFP allow us to construct
the ADF GoF tests as follows. We also discuss  the possibility of
construction of   similar tests in the case of i.i.d. observations and in
the case of nonlinear AR time series.

 First we show that the corresponding SFP's
$U_\varepsilon\left(\cdot \right) ,U_T\left(\cdot \right)$ and $U_n\left(\cdot
\right)$ converge to the  processes $ \left(U\left(\vartheta ,t\right), 0\leq t\leq
T\right)$, $ \left(U\left(\vartheta ,x\right), x\in R\right)$ and $
\left(U\left(\vartheta ,t\right),\right. $ $\left. 0\leq t\leq
\tau_*\right) $ respectively. Say, $U_\varepsilon \left(\cdot \right)$
converges to
$$
U\left(\vartheta ,t\right)=\int_{0}^{t}h\left(s\right){\rm
  d}W_s-\int_{0}^{T}h\left(s\right){\rm
  d}W_s\int_{0}^{t}h\left(s\right)^2{\rm d}s ,\;
\int_{0}^{T}h\left(s\right)^2{\rm d}s =1,
$$
where $h\left(s\right)=h\left(\vartheta ,s\right) $ is some function and
$W_s,0\leq s\leq T
$ is a Wiener process.
Therefore  if we put
$$
\tau =\int_{0}^{t}h\left(\vartheta ,s\right)^2{\rm d}s,\qquad \int_{0}^{t}h\left(s\right){\rm
  d}W_s=W\left(\int_{0}^{t}h\left(\vartheta ,s\right)^2{\rm d}s \right)=W\left(\tau \right),
$$
where $W\left(\cdot \right)$ is another Wiener process, then we can write
$$
U\left(\vartheta ,t\right)=W\left(\tau
\right)-W\left(1\right)\,\tau=B\left(\tau \right) ,\qquad 0\leq \tau \leq 1,
$$
where  $B\left(\cdot
\right)$ is a Brownian bridge. Hence
$$
\int_{0}^{T}U\left(\vartheta ,t\right)^2h\left(\vartheta ,t\right)^2{\rm
  d}t=\int_{0}^{1}B\left(\tau \right)^2{\rm d}\tau  =\Delta.
$$

This suggests the construction of  tests with the help of ``empirical
versions'' $U_{\varepsilon ,T,n}\left(\cdot \right)$ and $h_{\varepsilon
  ,T,n}\left(\cdot \right)$ of $U\left(\cdot \right)$ and $h\left(\cdot
\right)$ as
follows. Introduce the corresponding statistics (symbolic writing)
\begin{align*}
\Delta _{\varepsilon ,T,n}=\int_{}^{}U_{\varepsilon ,T,n}
\left(s\right)^2h_{\varepsilon ,T,n} \left(s\right)^2{\rm d}s .
\end{align*}
Then we show that for all three models we have the convergences to the same
limit
$$
\Delta _\varepsilon\Longrightarrow  \Delta ,\quad \Delta _T\Longrightarrow  \Delta
,\quad \Delta _n\Longrightarrow \Delta
$$
and therefore the tests
$$
\hat\psi _\varepsilon =\1_{\left\{\Delta _\varepsilon >c_\alpha
  \right\}},\qquad \hat\psi _T =\1_{\left\{\Delta _T>c_\alpha \right\}},\qquad
\hat\psi _n =\1_{\left\{\Delta _n >c_\alpha \right\}},\qquad
\Pb\left(\Delta >c_\alpha \right)=\alpha
$$
are ADF. Below we realize this program. Moreover we show that this approach
cannot be applied directly to the model of observations of i.i.d. random
variables, but in the case of nonlinear AR time series we have the similar ADF
GoF test, of course, under the strong regularity conditions.

This work is a continuation of the study of GoF tests for diffusion processes
observed in continuous time.  The case of simple basic hypothesis was treated
for example in the works \cite{Fou92},\cite{IK01}, \cite{Kut04}, \cite{DK},
\cite{NN}, \cite{Kut12b}.  The case of parametric basic hypothesis and ADF
tests was studied in the works \cite{NZ}, \cite{Kut12b}, \cite{KK},
\cite{Kut13a}, \cite{Kut13b}.

For point processes there are many publications devoted to this subject, see,
e.g., \cite{MTP} and the references therein.

\section{Score-Function Processes}

We have three stochastic processes observed in continuous time : small noise
diffusion, ergodic diffusion and inhomogeneous Poisson processes. First we consider
limits of the SFP's, separately for these models of observations.
Then we show how these limits can be used for  construction of
the ADF GoF tests.

\subsection{Small Noise Diffusion Processes.}

 We observe a realization
$
X^\varepsilon =\left(X_t,0\leq t\leq
  T\right)
$
of diffusion process  satisfying
  the stochastic differential equation
\begin{equation}
\label{1.1}
{\rm d}X_t=S\left(t,X_t\right){\rm d}t+\varepsilon \sigma \left(t,X_t\right)\,{\rm
  d}W_t, \quad x_0,\quad 0\leq t\leq T,
\end{equation}
where the trend coefficient $S\left(t,X_t\right)$ is an  unknown function and the
diffusion coefficient  $\varepsilon ^2\sigma \left(t,X_t\right)^2$  is a known
positive function. The initial value $x_0$ is deterministic and
$\varepsilon \in (0,1]$.

 We have to test the
following parametric (basic)  hypothesis:

\bigskip

${\scr H}_0$ : {\it The observed process
has the stochastic differential}
\begin{equation}
\label{1.2}
{\rm d}X_t=S\left(\vartheta ,t,X_t\right){\rm d}t+\varepsilon \sigma \left(t,X_t\right)\,{\rm
  d}W_t, \quad x_0,\quad 0\leq t\leq T,
\end{equation}
{\it where  the trend coefficient $S\left(\vartheta ,t,X_t\right)$ is a
  known smooth  function
depending on some unknown parameter $\vartheta \in \Theta =\left(a ,b \right)$. }

\bigskip

Our goal is to construct a GoF test $\hat\psi _\varepsilon $, which
belongs to the class ${\cal K}_\alpha $ and is consistent
in the asymptotics of {\it small noise }   $\varepsilon  \rightarrow 0 $. Note
that this stochastic model and the statistical inference for it has been considered in many works.
See, for example, \cite{FV}, \cite{Kut94} \cite{NY96} and the
references therein.

 Let us introduce the following regularity condition.

${\cal R}.$ {\it The functions $S\left(\vartheta ,t,x\right)$ and $\sigma
\left(t,x\right)$ have two continuous bounded derivatives with respect to
$\vartheta $ and $x$ and have continuous bounded derivatives w.r.t. $t$.}

 Below the
dot stands for the derivative w.r.t. $\vartheta $ and
prime means the derivative w.r.t. $x$ or w.r.t. $t$. For example,
$$
\ddot S\left(\vartheta ,t,x\right)=\frac{\partial^2 S\left(\vartheta ,t,x\right)}{\partial
\vartheta^2 },\quad S'_x\left(\vartheta ,t,x\right)=\frac{\partial
  S\left(\vartheta ,t,x\right)}{\partial
x} .
$$
Let us denote by  $x^T=\left(x_t,0\leq t\leq T\right)$ the solution of the
equation \eqref{1.2} with  $\varepsilon =0$, i.e. $x^T $  is solution of the   ordinary
differential equation
$$
\frac{{\rm d}x_t}{{\rm d}t}=S\left(\vartheta ,t,x_t\right),\qquad x_0,\quad 0 \leq t\leq T.
$$
Of course it is a function of $\vartheta$, i.e. $x_t=x_t\left(\vartheta \right)$.
It is known that as $\varepsilon \rightarrow 0$,  the process
$X^\varepsilon $ converges to the deterministic
function $x^T$
and this convergence is uniform w.r.t. $t\in \left[0,T\right]$ (see \cite{FV}).

Further, assume that the following identifiability condition is fulfilled.

${\cal I}.$  {\it For any} $\nu
>0$
$$
\inf_{\vartheta _0\in \Theta
  }\inf_{\left|\vartheta -\vartheta _0\right|>\nu
}\int_{0}^{T}\left[\frac{S\left(\vartheta ,t,x_t^*\right)-S\left(\vartheta_0,t
    ,x_t^*\right)}{\sigma \left(t,x_t^*\right)}\right]^2{\rm d}t >0.
$$
Here and below $x_t^*=x_t\left(\vartheta _0\right)$.

 The likelihood ratio function in the case of
observations \eqref{1.2} is
$$
L\left(\vartheta ,X^\varepsilon
\right)=\exp\left\{\int_{0}^{T}\frac{S\left(\vartheta
  ,t,X_t\right)}{\varepsilon ^2\sigma \left(t,X_t\right)^2}\;{\rm
  d}X_t-\int_{0}^{T}\frac{S\left(\vartheta
  ,t,X_t\right) ^2}{2\varepsilon ^2\sigma \left(t,X_t\right)^2}\;{\rm
  d}t\right\} ,\;\vartheta \in \Theta
$$
and the MLE  $\hat\vartheta _\varepsilon $ is defined by the equation
\begin{equation}
\label{1.3}
L\left(\hat\vartheta _\varepsilon ,X^\varepsilon \right)=\sup_{\vartheta \in
  \Theta }L\left(\vartheta ,X^\varepsilon \right).
\end{equation}

The MLE $\hat\vartheta _\varepsilon
$ under the aforementioned regularity conditions admits
the representation
\begin{equation}
\label{1.4}
\varepsilon ^{-1}\left(\hat\vartheta _\varepsilon -\vartheta
\right)={\rm I}\left(\vartheta \right)^{-1}\int_{0}^{T}\frac{\dot
  S\left(\vartheta ,t,x_t\right)}{\sigma \left(t,x_t\right)}\;{\rm d}W_t
+o\left(1\right)
\end{equation}
see \cite{Kut94}. Here   $ {\rm I}\left(\vartheta
\right)$  is the Fisher information
$$
{\rm I}\left(\vartheta
\right)=\int_{0}^{T}\left(\frac{\dot S\left(\vartheta ,t,x_t\right)}{\sigma
    \left(t,x_t\right)}\right)^2 \,{\rm d}t>0.
$$
We define the score-function
$$
\frac{\partial \ln L\left(\vartheta ,X^\varepsilon \right)}{\partial \vartheta
}=\int_{0}^{T}\frac{\dot S\left(\vartheta ,t,X_t\right)}{\varepsilon ^2\sigma
  \left(t,X_t\right)^2}\left[{\rm d}X_t-S\left(\vartheta ,t,X_t\right)\,{\rm
    d}t\right]
$$
and the normalized score-function
$$
U_\varepsilon\left(\vartheta ,X^\varepsilon \right)=\int_{0}^{T} \frac{\dot
  S\left(\vartheta
  ,t,X_t\right)}{\varepsilon \; {\rm
  I}\left(\vartheta \right)^{1/2}  \sigma \left(t,X_t\right)^2}\;\left[{\rm
    d}X_t-S\left(\vartheta ,t,X_t\right)\,{\rm d}t\right] .
$$
If the true value is $\vartheta _0$, then we have the convergence
$$
U_\varepsilon\left(\vartheta _0,X^\varepsilon \right)=\int_{0}^{T} \frac{\dot
  S\left(\vartheta_0
  ,t,X_t\right)}{ {\rm
  I}\left(\vartheta_0 \right)^{1/2}  \sigma \left(t,X_t\right)}\;{\rm
  d}W_t\longrightarrow  \zeta ,
$$
where
$$
\zeta  =\int_{0}^{T} \frac{\dot S\left(\vartheta_0
  ,t,x^*_t\right)}{ {\rm
  I}\left(\vartheta_0 \right)^{1/2}   \sigma \left(t,x^*_t\right)}\;{\rm
  d}W_t\quad \sim \quad {\cal N}\left(0,1 \right).
$$
The proof, which can be found in   \cite{Kut94}, follows from the uniform convergence of $X_t$ to $x_t^*$.

Let us introduce the score-function process
$$
U_\varepsilon\left(t,\vartheta ,X^\varepsilon \right)={\rm
  I}\left(\vartheta\right)^{-1/2}\int_{0}^{t} \frac{\dot S\left(\vartheta
  ,s,X_s\right)}{ \varepsilon \,   \sigma \left(s,X_s\right)^2}\;\left[{\rm
  d}X_s-S\left(\vartheta   ,s,X_s\right){\rm d}s\right],\quad 0\leq t\leq T,
$$
and (formally) the statistic
$
 U_\varepsilon \left(t\right)=U_\varepsilon\left(t,\hat\vartheta_\varepsilon
,X^\varepsilon \right),\quad 0\leq t\leq T.
$
We say ``formally'' because the MLE $\hat\vartheta _\varepsilon $ depends on
the whole trajectory $X^\varepsilon $ and the corresponding It\^o integral
\begin{equation}
\label{i}
\int_{0}^{t} \frac{\dot S(\hat\vartheta _\varepsilon
  ,s,X_s)}{  \sigma \left(s,X_s\right)^2}\;{\rm
  d}X_s
\end{equation}
is not well defined. The correct
definition will be given later and here we show (as well formally) to which  limit   this
process can be expected to converge.

 Note that $U_\varepsilon\left(T,\vartheta_0 ,X^\varepsilon
\right)=U_\varepsilon\left(\vartheta_0 ,X^\varepsilon \right)$ with
$\Pb_{\vartheta_0} ^{\left(\varepsilon \right)} $ probability 1.

  We have ($\vartheta _0$ is the true value)
\begin{align}
 U_\varepsilon \left(t\right)&=\int_{0}^{t} \frac{\dot S(\hat\vartheta _\varepsilon
  ,s,X_s)}{ \varepsilon \, {\rm
  I}(\hat\vartheta _\varepsilon )^{1/2}\;  \sigma \left(s,X_s\right)^2}\;\left[{\rm
  d}X_s-S(\hat\vartheta _\varepsilon
  ,s,X_s)\;{\rm d}s\right]\nonumber\\
&=\int_{0}^{t} \frac{\dot S(\hat\vartheta _\varepsilon
  ,s,X_s)}{  {\rm
  I}(\hat\vartheta _\varepsilon )^{1/2}\;  \sigma \left(s,X_s\right)}\;{\rm
  d}W_s\nonumber\\
&\qquad +\int_{0}^{t} \frac{\dot S(\hat\vartheta _\varepsilon
  ,s,X_s)\left[S(\vartheta _0
  ,s,X_s)-S(\hat\vartheta _\varepsilon
  ,s,X_s)\right]}{ \varepsilon \, {\rm
  I}(\hat\vartheta _\varepsilon )^{1/2}\;  \sigma \left(s,X_s\right)^2}\;{\rm d}s\nonumber\\
&=\int_{0}^{t} \frac{\dot S(\hat\vartheta _\varepsilon
  ,s,X_s)}{  {\rm
  I}(\hat\vartheta _\varepsilon  )^{1/2}  \sigma \left(s,X_s\right)}\;{\rm
  d}W_s\nonumber\\
&\qquad -\frac{\hat\vartheta _\varepsilon -\vartheta _0  }{\varepsilon\; {\rm
  I}(\hat\vartheta _\varepsilon  )^{-1/2} }\int_{0}^{t} \frac{\dot S(\hat
  \vartheta _\varepsilon
  ,s,X_s)\dot S(\tilde \vartheta _\varepsilon
  ,s,X_s)}{  {\rm
  I}(\hat\vartheta _\varepsilon ) \; \sigma \left(s,X_s\right)^2}\;{\rm d}s\nonumber\\
&=\int_{0}^{t} \frac{\dot S(\vartheta _0
  ,s,x^*_s)}{  {\rm
  I}(\vartheta _0 )^{1/2}\;  \sigma \left(s,x^*_s\right)^2}\;{\rm
  d}W_s\nonumber\\
&\qquad -\int_{0}^{T} \frac{\dot S(\vartheta _0
  ,s,x^*_s)}{  {\rm
  I}(\vartheta _0 )^{1/2}  \sigma \left(s,x^*_s\right)}\;{\rm
  d}W_s\int_{0}^{t} \frac{\dot S(\vartheta _0
  ,s,x^*_s)^2}{  {\rm
  I}(\vartheta _0 ) \; \sigma \left(s,x^*_s\right)^2}\;{\rm
  d}s+o\left(1\right).
\label{1.7}
\end{align}
Further, if we denote
$$
\tau =\int_{0}^{t} \frac{\dot S(\vartheta _0
  ,s,x^*_s)^2}{  {\rm
  I}(\vartheta _0 )  \sigma \left(s,x^*_s\right)^2}\;{\rm d}s,\qquad 0\leq \tau \leq 1,
$$
 then we can write
$$
\int_{0}^{t} \frac{\dot S(\vartheta _0
  ,s,x^*_s)}{  {\rm
  I}(\vartheta _0 )^{1/2} \; \sigma \left(s,x^*_s\right)^2}\;{\rm
  d}W_s=W\left(\tau \right),
$$
where $W\left(\cdot \right)$ is some Wiener process.  Therefore we obtain the
limit
$$
 U_\varepsilon \left(t\right)\Longrightarrow W\left(\tau \right)-W\left(1
\right)\tau =B\left(\tau \right),\qquad 0\leq \tau \leq 1,
$$
with a Brownian bridge $B\left(\cdot \right)$.

This convergence suggests the construction of the following test statistic
\begin{equation}
\label{qu}
\Delta _\varepsilon =\int_{0}^{T}\frac{ U_\varepsilon \left(t\right)^2\dot
  S(\hat\vartheta _\varepsilon    ,t,X_t)^2}{  {\rm
  I}(\hat\vartheta _\varepsilon  )\;  \sigma \left(t,X_t\right)^2}\;{\rm d}t
\end{equation}
and the test
$
\hat\psi _\varepsilon =\1_{\left\{\Delta _\varepsilon >c_\alpha
  \right\}}$,where   $ \Pb\left(\Delta >c_\alpha \right)=\alpha .
$
If we verify that
$$
\Delta _\varepsilon\Longrightarrow \Delta =\int_{0}^{1}B\left(\tau
\right)^2{\rm d}\tau ,
$$
then the test $\hat\psi _\varepsilon\in {\cal K}_\alpha  $ and is ADF.

To avoid the problem concerning  the stochastic integral \eqref{i} we use two
possibilities: one is the
well-known device which consists
in the application of the It\^o formula to the function
$$
H\left(\vartheta ,s,x\right)=\int_{x_0}^{x}\frac{\dot S\left(\vartheta
  ,s,y\right)}{\sigma \left(s,y\right)^2}{\rm d}y
$$
and the second is based on some preliminary estimator of the parameter
$\vartheta $.

The first approach was applied in the similar problem  in \cite{Kut13a} and
here we follow the same steps. The second approach was mentioned in
\cite{Kut13a} too but here (below) we work out  the details of the proof.

{\it The first approach.}
The It\^o formula  applied to  the function $H\left(\vartheta ,s,X_s\right)$
gives us the stochastic differential
\begin{align*}
\int_{0}^{t}\frac{\dot S\left(\vartheta ,s,X_s\right)}{\sigma
  \left(s,X_s\right)^2}\,{\rm d}X_s&=
H\left(\vartheta ,t,X_t\right)\\
&\quad -\int_{0}^{t}\left[H'_s\left(\vartheta
  ,s,X_s\right)+\frac{\varepsilon^2\,\sigma \left(s,X_s\right)^2}{2
  }H''_{x,x}\left(\vartheta ,s,X_s\right)\right]{\rm d}s.
\end{align*}
Note that the contribution of the term
$$
\varepsilon^2\,\int_{0}^{t}{\sigma \left(s,X_s\right)^2}{
  }H''_{x,x}\left(\vartheta ,s,X_s\right){\rm d}s
$$
is asymptotically  negligible and we can omit it.

We have
\begin{align}
\hat U_\varepsilon\left(t \right)&= \frac{ H(\hat\vartheta _\varepsilon ,t,X_t)}{\varepsilon {\rm
I}(\hat\vartheta _\varepsilon )^{1/2} } -\int_{0}^{t}\frac{H'_s\left(\hat\vartheta _\varepsilon
  ,s,X_s\right)}{\varepsilon{\rm
I}(\hat\vartheta _\varepsilon )^{1/2}  }{\rm d}s-\int_{0}^{t} \frac{\dot
  S(\hat\vartheta _\varepsilon
  ,s,X_s)S(\vartheta _0
  ,s,X_s)}{ \varepsilon  {\rm
  I}(\hat\vartheta _\varepsilon )^{1/2}\;  \sigma \left(s,X_s\right)^2}{\rm d}s\nonumber\\
&\qquad -\int_{0}^{t} \frac{\dot S(\hat\vartheta _\varepsilon
  ,s,X_s)\left[S(\hat\vartheta _\varepsilon
  ,s,X_s)-S(\vartheta _ 0
  ,s,X_s)\right]}{ \varepsilon \, {\rm
  I}(\hat\vartheta _\varepsilon )^{1/2}\;  \sigma \left(s,X_s\right)^2}\;{\rm
  d}s+O\left(\varepsilon \right)\nonumber\\
&=J\left(\hat\vartheta _\varepsilon,t,X^t  \right)-K\left(\hat\vartheta
_\varepsilon,t,X^t  \right)+O\left(\varepsilon \right),
\label{u}
\end{align}
where $K\left(\cdot \right)$ is the last integral.  Its convergence is obtained directly (see \eqref{1.4}):
\begin{align*}
K\left(\hat\vartheta _\varepsilon,t,X^t  \right)&=\frac{\hat\vartheta
  _\varepsilon-\vartheta _0 }{\varepsilon }\; \int_{0}^{t} \frac{\dot
  S(\hat\vartheta _\varepsilon
  ,s,X_s)\dot S(\tilde\vartheta _\varepsilon
  ,s,X_s)}{  {\rm
  I}(\hat\vartheta _\varepsilon )^{1/2}\;  \sigma \left(s,X_s\right)^2}\;{\rm
  d}s\\
&\longrightarrow   \int_{0}^{T} \frac{\dot S(\vartheta _0
  ,s,x^*_s)}{ {\rm
  I}(\vartheta _0  )^{1/2}\; \sigma \left(s,x^*_s\right)}\;{\rm
  d}W_s\;\int_{0}^{t} \frac{\dot S(\vartheta _0
  ,s,x^*_s)^2}{  {\rm
  I}(\vartheta _0  )\;  \sigma \left(s,x^*_s\right)^2}\;{\rm   d}s.
\end{align*}
Further, we verify that
\begin{align*}
J\left(\hat\vartheta _\varepsilon,t,X^t  \right)-J\left(\vartheta _0,t,X^t
\right)\longrightarrow 0
\end{align*}
and that
$$
J\left(\vartheta _0,t,X^t  \right)\longrightarrow \int_{0}^{t} \frac{\dot
  S(\vartheta _0
  ,s,x^*_s)}{ {\rm
  I}(\vartheta _0  )^{1/2}\; \sigma \left(s,x^*_s\right)}\;{\rm
  d}W_s
$$
(see details in \cite{Kut13a}).

Thus we obtained the convergence mentioned in \eqref{1.7} and the
following result.
\begin{proposition}
\label{P1} Suppose that the conditions of regularity are fulfilled, then the
test $\hat\psi _\varepsilon =\1_{\left\{\Delta _\varepsilon >c_\alpha
  \right\}}$ with
\begin{align*}
\Delta _\varepsilon =\int_{0}^{T}\frac{\hat U_\varepsilon \left(t\right)^2\dot
  S(\hat\vartheta _\varepsilon    ,t,X_t)^2}{  {\rm
  I}(\hat\vartheta _\varepsilon  )\;  \sigma \left(t,X_t\right)^2}\;{\rm d}t
\end{align*}
is ADF and belongs to ${\cal K}_\alpha $.
\end{proposition}

{\it Second approach.}
Let us write $\hat U_\varepsilon \left(t\right)$ as the difference  of two integrals
\begin{align*}
\hat U_\varepsilon \left(t\right)=\int_{0}^{t} \frac{\dot S(\hat\vartheta _\varepsilon
  ,s,X_s)}{\varepsilon\, {\rm I}(\hat\vartheta
  _\varepsilon )^{1/2} \sigma \left(s,X_s\right)^2}\;{\rm
  d}X_s-\int_{0}^{t} \frac{\dot S(\hat\vartheta _\varepsilon
  ,s,X_s)\;S(\hat\vartheta _\varepsilon
  ,s,X_s)}{  \varepsilon\, {\rm I}(\hat\vartheta _\varepsilon )^{1/2} \sigma
  \left(s,X_s\right)^2}\;{\rm d}s.
\end{align*}
Note that the properties of the estimator $\hat\vartheta _\varepsilon $
required  in the study of the
first and the second integrals are different.

In the first integral it is sufficient that
$\hat\vartheta _\varepsilon  \rightarrow \vartheta_0 $ and in the second
integral we need the asymptotic efficiency (full limit variance) of the
MLE. Therefore we can consider two different estimators in the calculation of these
integrals. For the first integral we introduce a preliminary (consistent)
estimator $\bar\vartheta _{\nu _\varepsilon }$ constructed by the first $\left(X_t,
0\leq t\leq \nu _\varepsilon\right) $ observations. Here $\nu _\varepsilon
\rightarrow 0$ but slowely. Then  we can
use the estimator $\bar\vartheta _{\nu _\varepsilon } $ in the calculation of the integral
$$
\int_{\nu _\varepsilon }^{t} \frac{\dot S(\bar\vartheta _{\nu _\varepsilon}
  ,s,X_s)}{ \varepsilon\, {\rm I}(\bar\vartheta _{\nu _\varepsilon}
  )^{1/2}\sigma \left(s,X_s\right)^2}\;{\rm
  d}X_s,\qquad t\in \left[\nu _\varepsilon ,T\right],
$$
which is now well defined. In the second integral we keep $\hat\vartheta
_\varepsilon $ in the function $S(\hat\vartheta _\varepsilon
  ,s,X_s)$ only. Therefore we  consider the statistic
\begin{align*}
V_\varepsilon \left(t\right)=\int_{\nu _\varepsilon}^{t} \frac{\dot S(\bar\vartheta _{\nu _\varepsilon}
  ,s,X_s)}{\varepsilon\, {\rm I}(\bar\vartheta
  _{\nu _\varepsilon} )^{1/2} \sigma \left(s,X_s\right)^2}\;{\rm
  d}X_s-\int_{\nu _\varepsilon}^{t} \frac{\dot S(\bar\vartheta
  _{\nu _\varepsilon}
  ,s,X_s)\;S(\hat\vartheta _\varepsilon
  ,s,X_s)}{  \varepsilon\, {\rm I}(\bar\vartheta
  _{\nu _\varepsilon} )^{1/2} \sigma
  \left(s,X_s\right)^2}\;{\rm d}s,
\end{align*}
where $ t\in \left[\nu _\varepsilon ,T\right]$. Now we can repeat the
calculations similar to  \eqref{1.7} for the statistic $V_\varepsilon
\left(t\right),t\in \left[\nu _\varepsilon ,T\right] $, which is this time well defined, and obtain the same limit
expression.

Let us construct a consistent estimator $\bar\vartheta _{\nu _\varepsilon} $
by the ``vanishing observations''  $ X_t,
0\leq t\leq \nu _\varepsilon, \nu _\varepsilon \rightarrow 0  $. Introduce a
minimum distance estimator (MDE)
\begin{align*}
\bar\vartheta _{\nu _\varepsilon }=\arg \inf_{\vartheta \in \Theta }\int_{0}^{\nu _\varepsilon }\left[X_t-x_t\left(\vartheta \right)\right]^2{\rm d}t.
\end{align*}
 The consistency of this estimator is verified
in the following lemma.

\begin{lemma}
\label{L1} Suppose that the regularity condition ${\cal R}$ is fulfilled and
for all $\vartheta \in \Theta $ we have
$\left|\dot
S\left(\vartheta ,0,x_0\right)\right|\geq \kappa $, where $\kappa >0$. Then
the MDE $\bar\vartheta _{\nu _\varepsilon} $ with $\nu _\varepsilon
=\varepsilon^2 \ln \left(\varepsilon ^{-1}\right)$ is consistent.
\end{lemma}
{\bf Proof.} Below $\left\|\cdot \right\|_{\nu _\varepsilon }$ is $L^2\left[0,\nu
  _\varepsilon \right]$ norm. Let us put
$$
g\left(\gamma ,\nu _\varepsilon \right)=\inf_{\left|\vartheta -\vartheta
  _0\right|>\gamma }\left\|x_t\left(\vartheta \right)-x_t\left(\vartheta_0
\right)\right\|_{\nu _\varepsilon } .
$$
 Note that
\begin{align*}
g\left(\gamma ,\nu _\varepsilon \right) ^2 &=\int_{0}^{\nu _\varepsilon }\left[x_t\left(\vartheta
  \right)-x_t\left(\vartheta_0 \right)\right] ^2{\rm d}t=\left(\vartheta
-\vartheta _0\right)^2 \int_{0}^{\nu _\varepsilon }\dot x_t(\tilde \vartheta
  )^2{\rm d}t.
\end{align*}
with some $\tilde \vartheta  $. The
derivative w.r.t. $\vartheta $ of $x_t\left(\vartheta \right)$ satisfies the
equation
\begin{align*}
\frac{{\rm d}\dot x_t\left(\vartheta \right)}{{\rm d}t}=\dot S\left(\vartheta,t
,x_t\left(\vartheta \right)\right) +S_x'\left(\vartheta ,t,x_t\left(\vartheta
\right)\right)\dot x_t\left(\vartheta \right),\quad \dot x_0\left(\vartheta \right)=0.
\end{align*}
Its solution is the function
\begin{align*}
\dot x_t\left(\vartheta \right)=\int_{0}^{t}\dot S\left(\vartheta,s
,x_s\left(\vartheta \right)\right) \exp\left\{\int_{s}^{t}S'_x\left(\vartheta
,v,x_v\left(\vartheta \right)\right){\rm d}v\right\} \;{\rm d}s.
\end{align*}
Hence for the small values of $t$ we have the estimate
\begin{align*}
\dot x_t\left(\vartheta \right)=t\dot S\left(\vartheta,0
,x_0\right)\left(1+O\left(t\right)\right).
\end{align*}
Therefore for all $\varepsilon <\varepsilon _*$, where
$\varepsilon _*$ is some small value
\begin{align*}
\left\|x_t\left(\vartheta \right)-x_t\left(\vartheta_0
\right)\right\|_{\nu _\varepsilon } ^2 \geq \frac{\left(\vartheta
-\vartheta _0\right)^2\kappa ^2\nu _\varepsilon ^3}{6}.
\end{align*}

Further, for any $\gamma >0$ we have
\begin{align*}
&\Pb_{\vartheta _0}\left(\left|\bar\vartheta _{\nu _\varepsilon }-\vartheta
_0\right|>\gamma \right)\\
&\qquad \quad =\Pb_{\vartheta _0} \left(\inf_{\left|\vartheta
 -\vartheta _0\right|\leq \gamma } \left\|X_t-x_t\left(\vartheta
\right)\right\|_{\nu _\varepsilon }      > \inf_{\left|\vartheta
 -\vartheta _0\right|> \gamma } \left\|X_t-x_t\left(\vartheta
\right)\right\|_{\nu _\varepsilon }  \right)\\
&\qquad \quad  \leq  \Pb_{\vartheta _0} \left(\inf_{\left|\vartheta
 -\vartheta _0\right|\leq \gamma }\left( \left\|X_t-x_t\left(\vartheta_0
\right)\right\|_{\nu _\varepsilon }+\left\|x_t\left(\vartheta \right)-x_t\left(\vartheta_0
\right)\right\|_{\nu _\varepsilon } \right)  \right.\\
&\qquad  \qquad \qquad\left. > \inf_{\left|\vartheta
 -\vartheta _0\right|> \gamma } \left(\left\|x_t\left(\vartheta \right)-x_t\left(\vartheta_0
\right)\right\|_{\nu _\varepsilon }-\left\|X_t-x_t\left(\vartheta_0
\right)\right\|_{\nu _\varepsilon } \right) \right)\\
&\qquad \quad =\Pb_{\vartheta _0} \left( 2\left\|X_t-x_t\left(\vartheta_0
\right)\right\|_{\nu _\varepsilon } \geq g\left(\gamma ,\nu _\varepsilon
\right)\right)\\
&\qquad \qquad \leq \frac{4}{g\left(\gamma ,\nu _\varepsilon
\right)^2} \Ex_{\vartheta _0}\int_{0}^{\nu _\varepsilon
}\left[X_t-x_t\left(\vartheta _0\right)\right]^2{\rm d}t \leq
\frac{C\varepsilon ^2\nu _\varepsilon^2 }{\gamma ^2\kappa ^2\nu _\varepsilon ^3 }\leq \frac{C}{\ln
 \frac{1}{ \varepsilon} }\longrightarrow 0.
\end{align*}
Here we used the estimate
\begin{align*}
\sup_{0\leq s\leq t}\Ex_{\vartheta _0} \left|X_s-x_s\left(\vartheta _0 \right)\right|^2\leq Ct\varepsilon ^2,
\end{align*}
which can be found, for example, in \cite{Kut94}, Lemma 1.13.

Therefore the estimator $\bar\vartheta _{\nu _\varepsilon} $ is consistent and
we have the following result.
\begin{proposition}
\label{P2} Suppose that the conditions of regularity are fulfilled and for all
$\vartheta \in \Theta $ we have
$\left|\dot
S\left(\vartheta ,0,x_0\right)\right|\geq \kappa $, where $\kappa >0$, then the
test $\tilde\psi _\varepsilon =\1_{\left\{\tilde\Delta _\varepsilon >c_\alpha
  \right\}}$ with
\begin{align*}
\tilde\Delta _\varepsilon=\int_{\nu _\varepsilon }^{T}\frac{V_\varepsilon
  \left(t\right)^2\dot S\left(\bar\vartheta _{\nu _\varepsilon
  },t,X_t\right)^2}{{\rm I}\left(\bar\vartheta _{\nu _\varepsilon
  }  \right)\sigma \left(t,X_t\right)^2}\,{\rm d}t
\end{align*}
is ADF and belongs to ${\cal K}_\alpha $.
\end{proposition}

\bigskip

Let us consider the problem of  {\it  consisteny} of this test.  The observed
process under alternative is
$$
{\rm d}X_t=S\left(t,X_t\right){\rm d}t+\varepsilon \sigma \left(t,X_t\right){\rm
  d}W_t,\quad X_0=x_0,\quad 0\leq t\leq T ,
$$
where $S\left(t,x\right)$ does not belong to the parametric family of trend
coefficients $ \left\{S\left(\vartheta ,t,x\right),\vartheta \in \Theta \right\}$. We obtain
the following representation for the statistic
$V_\varepsilon \left(\cdot \right)$:
\begin{align*}
 V_\varepsilon \left(t\right)&=\int_{\nu _\varepsilon}^{t} \frac{\dot S(\bar\vartheta _{\nu _\varepsilon}
  ,s,X_s)}{  {\rm
  I}(\bar\vartheta _{\nu _\varepsilon} )^{1/2}\;  \sigma \left(s,X_s\right)^2}\;{\rm
  d}W_s\\
&\qquad +\int_{{\nu _\varepsilon}}^{t} \frac{\dot S(\bar\vartheta _{\nu _\varepsilon}
  ,s,X_s)\left[S(s,X_s)-S(\hat\vartheta _\varepsilon
  ,s,X_s)\right]}{ \varepsilon \, {\rm
  I}(\bar\vartheta _{\nu _\varepsilon} )^{1/2}\;  \sigma \left(s,X_s\right)^2}\;{\rm
  d}s\\
&=\int_{{\nu _\varepsilon}}^{t} \frac{\dot S(\bar\vartheta
  ,s,x_s)}{  {\rm
  I}(\bar\vartheta )^{1/2}\;  \sigma \left(s,x_s\right)^2}\;{\rm
  d}W_s+o\left(1\right)\\
&\qquad +\int_{{\nu _\varepsilon}}^{t} \frac{\dot S(\bar\vartheta
  ,s,x_s)\left[S(s,x_s)-S(\hat\vartheta
  ,s,x_s)\right]}{ \varepsilon \, {\rm
  I}(\bar\vartheta)^{1/2}\;  \sigma \left(s,x_s\right)^2}\;{\rm
  d}s\left(1+o\left(1\right)\right).
\end{align*}
Here $x_t$ is solution of the ordinary differential  equation
$$
\frac{{\rm d}x_t}{{\rm d}t}=S\left(t,x_t\right),\qquad x_0,\quad 0\leq t\leq T
$$
and $\hat\vartheta,\bar\vartheta $ are  defined as follows
\begin{align}
\label{misml}
\hat\vartheta&=\arg \inf_{\vartheta\in \Theta } \int_{0}^{T}\left(\frac{S\left(\vartheta
  ,t,x_t\right)-S\left(t,x_t\right)}{\sigma \left(t,x_t\right)}\right)^2{\rm d}t,\\
\bar\vartheta&=\arg \inf_{\vartheta\in \Theta } \left|S\left(\vartheta
  ,0,x_0\right)-S\left(0,x_0\right) \right|.
\label{mismd}
\end{align}
For the proof of \eqref{misml} see \cite{Kut94}, Section 2.6 and the equality
\eqref{mismd} is obtained as follows.  We have
\begin{align*}
\left\|x_t-x_t\left(\vartheta \right)\right\|^2_{\nu _\varepsilon
}&=\int_{0}^{\nu _\varepsilon } \left[x_t-x_t\left(\vartheta
  \right)\right]^2{\rm d}t\\
&=\int_{0}^{\nu _\varepsilon } t^2\left[S\left(0,x_0\right)-S\left(\vartheta
  ,0,x_0\right)\right]^2{\rm d}t\,\left(1+o\left(1\right)\right).
\end{align*}
Hence
\begin{align*}
\bar\vartheta _{\nu _\varepsilon}&=\arg\inf_{\vartheta \in\Theta
}\left\|x_t-x_t\left(\vartheta \right)\right\|^2_{\nu _\varepsilon} \\
&=\arg\inf_{\vartheta \in\Theta
}\frac{\nu _\varepsilon ^3}{3} \left[S\left(0,x_0\right)-S\left(\vartheta
  ,0,x_0\right)\right]^2\left(1+o\left(1\right)\right)\longrightarrow \bar\vartheta ,
\end{align*}
which yields \eqref{mismd}.

Introduce  the condition
\begin{align*}
\inf_{\bar\vartheta ,\hat\vartheta }\sup_{0\leq t\leq T}\left|\int_{0}^{t}
\frac{\dot S\left(\bar\vartheta ,s,x_s\right)
  \left[S\left(s,x_s\right)-S\left(\hat\vartheta ,s,x_s\right) \right]
}{\sigma \left(s,x_s\right)^2}{\rm d}s \right|>0.
\end{align*}
 It is easy to see that if this condition is fulfilled then $\Delta
 _\varepsilon \rightarrow \infty $ and the test is consistent.
Note that if this condition is not fulfilled then  for all $t\in \left[0,T\right]$ we have
\begin{align*}
&\int_{0}^{t}  \frac{\dot S\left(\bar\vartheta
  ,s,x_s\right) \left[S\left(s,x_s\right)-S\left(\hat\vartheta
  ,s,x_s\right)  \right] }{\sigma \left(s,x_s\right)^2}{\rm d}s=0
\end{align*}
and this equality implies
\begin{equation}
\label{as}
\dot S\left(\bar\vartheta
  ,t,x_t\right) \left[S\left(t,x_t\right)-S\left(\hat\vartheta
  ,t,x_t\right)  \right]=0,\qquad 0\leq t\leq T.
\end{equation}
If $\left|\dot S\left(\vartheta ,t,x\right)\right|>0$ for all $\vartheta \in
\Theta $ and almost all $t\in\left[0,T\right]$ and almost all $x\in K$ for any bounded
region $K\subset {\cal R}$, then the proposed test is consistent against any
fixed alternative.

An example of alternative invisible by this test can be constructed as
follows. Suppose that the function $S\left(\vartheta ,t,x\right)$ does not
depend on $\vartheta $ for the values $t\in \left[0,T/2\right]$ and the
trend coefficient $S\left(t,x_t\right)$  under alternative coincides with the
function $S\left(\vartheta^* ,t,x_t\right)$ for $t\in\left[T/2,
  T\right]$. Then we have \eqref{as} in the situation, where the trend
coefficients of diffusion process   on the interval $\left[0,T/2\right]$ can
be different under alternative. Of course as we know that the trend coefficient under hypothesis
does not depend on $\vartheta $ on the interval $\left[0,T/2\right]$, then for
this interval we can modify the test statistic.

{\bf Example.} Suppose that the observed diffusion process under hypothesis
has the stochastic differential
$$
{\rm d}X_t=\vartheta X_t\,{\rm d}t+\varepsilon\, {\rm d}W_t,\quad X_0=x_0>0,\quad 0\leq t\leq T,
$$
where $\vartheta \in \Theta $ and $0\not\in \Theta  $.
Then we have
$$
{\rm I}\left(\vartheta
\right)=\frac{x_0^2(e^{2\vartheta T }-1)}{2\vartheta},\qquad \hat\vartheta
_\varepsilon =\frac{\int_{0}^{T} X_t\,{\rm d}X_t}{\int_{0}^{T} X_t^2\,{\rm d}t}
 $$
and the statistic
\begin{align*}
\hat U_\varepsilon \left(t\right)=\frac{1}{\varepsilon x_0} {\sqrt{\frac{2\hat\vartheta _\varepsilon
    T}{ e^{2\hat\vartheta _\varepsilon T }-1   }}}
\int_{0}^{t}X_s\,\left[{\rm d}X_s-\hat\vartheta _\varepsilon X_t\,{\rm
    d}t\right].
\end{align*}
Here we have no problem of the definition of stochastic integral and this will always be  the case
for the models in which the trend coefficient depends linearly on the unknown parameter.

The test $\hat\psi _\varepsilon =\1_{\left\{\Delta _\varepsilon >c_\alpha
  \right\}}$ with
\begin{align*}
\Delta _\varepsilon =\int_{0}^{T}\frac{\hat U_\varepsilon \left(t\right)\;X_t^2
}{{\rm I}(\hat\vartheta_\varepsilon )\; \sigma ^2}\;{\rm d}t\Longrightarrow
\int_{0}^{1}B\left(\tau \right)^2\;{\rm d}\tau
\end{align*}
is ADF.

\bigskip

\subsection{Ergodic Diffusion Processes}

Suppose that the observed  diffusion process
$X^T=\left(X_t,0\leq t\leq T\right)$ satisfies the stochastic differential
\begin{equation}
\label{2.1}
{\rm d}X_t=S\left(X_t\right)\,{\rm d}t+ \sigma
\left(X_t\right)\,{\rm d}W_t, \quad X_0,\quad 0\leq t\leq T,
\end{equation}
where the function  $\sigma \left(x\right)$ is known. The trend
coefficient $S\left(\cdot \right)$ is an unknown  function and we have to test
the following composite   hypothesis:

\bigskip

${\scr H}_0$ : {\it The process $X^T$ is the solution of  equation}
\begin{equation}
\label{2.2}
{\rm d}X_t=S\left(\vartheta ,X_t\right)\,{\rm d}t+ \sigma
\left(X_t\right)\,{\rm d}W_t, \quad X_0,\quad 0\leq t\leq T,\quad
\vartheta \in \Theta ,
\end{equation}
{\it where $S\left(\vartheta ,x\right)$ is a known  smooth function depending on
unknown parameter $\vartheta \in \Theta =\left(a,b\right)$.}

\bigskip

Introduce the regularity conditions.

${\cal ES}.$ {\it The function $S\left(\vartheta ,x \right)$ is locally bounded, the
function $\sigma \left(\cdot \right)^2>0 $ is continuous  and for some $C>0$
the condition
$$
x\,S\left(\vartheta
  ,x\right)+\sigma \left(x\right)^2\leq C\left(1+x^2\right)
$$
holds.}

By this condition the stochastic differential equation has a unique weak
solution (see, e.g., \cite{Durett}).

Let us denote by  ${\cal P}$  the class of locally bounded functions with polynomial
majorants ($p>0$)
$$
{\cal P}=\left\{h\left(\cdot \right):\quad \left|h\left(y\right)\right|\leq
C\left(1+\left|y\right|^p\right)\right\} .
$$
  The next  condition  is

${\cal A}_0.$ {\it The functions  $S\left(\cdot \right), \sigma \left(\cdot
  \right)^{\pm 1} \in {\cal     P} $ and}
$$
\Limsup_{\left|y\right|\rightarrow \infty
}\;\sup_{\vartheta \in \Theta }\;\;\sgn\left(y\right)\;\frac{S\left(\vartheta
  ,y\right)}{\sigma \left(y\right)^2} <0.
$$
Note that if $S\left(\vartheta ,x \right)$ and $\sigma \left(x \right)$  satisfy
${\cal A}_0$, then we have
$$
V\left(\vartheta
,x\right)=\int_{0}^{x}\exp\left\{-2\int_{0}^{y}\frac{S\left(\vartheta
  ,z\right)}{\sigma   \left(z\right)^2}{\rm d}z\right\} {\rm
  d}y\longrightarrow \pm \infty
 $$
as $x\rightarrow \pm \infty $  and $\sup_{\vartheta \in \Theta
}G\left(\vartheta \right)<\infty  $, where
$$
G\left(\vartheta \right)=\int_{-\infty }^{\infty }\sigma
\left(y\right)^{-2}\exp\left\{2\int_{0}^{x}\frac{S\left(\vartheta ,y\right)}{\sigma
  \left(y\right)^2}{\rm d}y\right\} {\rm d}x
$$
is normalizing constant.

By these  conditions the stochastic process $X^T$
is positive-recurrent  (ergodic) with
the density of the invariant law
$$
f\left(\vartheta ,x\right)=\frac{1}{G\left(\vartheta \right)\;\sigma
  \left(x\right)^2}\;\exp\left\{2\int_{0}^{x}\frac{S\left(\vartheta ,y\right)}{\sigma
  \left(y\right)^2}\; {\rm d}y\right\}.
$$
Let us introduce further regularity conditions.

${\cal R}_e$.
{\it  The function $S\left(\vartheta ,x\right)$ has two continuous
derivatives }
$$
\dot S\left(\vartheta ,x\right),\ddot S\left(\vartheta
,x\right)\in {\cal P}.
$$
and

${\cal I}_e$.  {\it For any $\nu >0$}
$$
\inf_{\vartheta _0\in \Theta }\inf_{\left|\vartheta
  -\vartheta _0\right|>\nu } \int_{-\infty }^{\infty
}\left[\frac{S\left(\vartheta ,x\right)-S\left(\vartheta_0 ,x\right)}{\sigma
    \left(x\right)}\right]^2 f\left(\vartheta_0 ,x\right)\,{\rm d}x>0.
$$
 The likelihood ratio function is
\begin{align*}
L\left(\vartheta ,X^T\right)=\exp\left\{\int_{0}^{T}\frac{S\left(\vartheta
  ,X_t\right)}{\sigma \left(X_t\right)^2}\;{\rm d}X_t-\int_{0}^{T}\frac{S\left(\vartheta
  ,X_t\right)^2}{2\,\sigma \left(X_t\right)^2}\;{\rm d}t\right\} .
\end{align*}
Under the regularity conditions assumed above,
 the MLE $\hat\vartheta _T$  admits
 the representation
$$
\sqrt{T}\left(\hat\vartheta _T-\vartheta \right)=\frac{1}{{\rm I}\left(\vartheta
\right)\sqrt{T}} \int_{0}^{T}\frac{\dot S\left(\vartheta
  ,X_t\right)}{\sigma \left(X_t\right)} \; {\rm d}W_t+o\left(1\right).
$$
Here   $ {\rm I}\left(\vartheta
\right)$  is the Fisher information
$$
{\rm I}\left(\vartheta
\right)=\int_{-\infty }^{\infty
}\left(\frac{\dot S\left(\vartheta ,x\right)}{\sigma
    \left(x\right)}\right)^2 f\left(\vartheta ,x\right)\,{\rm d}x>0.
$$
The proof can be found in \cite{Kut04}.

The score-function is
\begin{align*}
\frac{\partial \ln L\left(\vartheta ,X^T\right)}{\partial \vartheta
}=\int_{0}^{T}\frac{\dot S\left(\vartheta
  ,X_t\right)}{\sigma \left(X_t\right)^2}\;\left[{\rm d}X_t-S\left(\vartheta
  ,X_t\right)\,{\rm d}t\right]
\end{align*}
and we define the normalized score-function:
$$
U_T\left(\vartheta,X^T\right)=\varphi _T\left(\vartheta
  \right)\int_{0}^{T}\frac{\dot S\left(\vartheta ,X_t\right)}{\sigma
  \left(X_t\right)^2}\;\left[{\rm d}X_t-S\left(\vartheta ,X_t\right)\,{\rm
      d}t\right] \Longrightarrow \xi ,
$$
where $\varphi _T\left(\vartheta
 \right)=\left[T{\rm I}\left(\vartheta
\right) \right]^{-1/2}$.
The limit random variable $\xi $ can be written as the following integral
$$
\xi =\int_{-\infty }^{\infty }\frac{\dot S\left(\vartheta
  ,y\right)\,\sqrt{f\left(\vartheta ,y\right)}}{\sqrt{{\rm
      I}\left(\vartheta \right)}\,\sigma
  \left(y\right)}\;{\rm d}w\left(y\right)\quad \sim\quad {\cal N}\left(0,1\right),
$$
where $w\left(\cdot \right)$ is two-sided Wiener process.

Let us introduce the slightly modified  score-function process
$$
U_T\left(x,\vartheta ,X^T\right)=\varphi _T\left(\vartheta
 \right)\int_{0}^{T}\frac{\dot S\left(\vartheta ,X_t\right)}{\sigma
  \left(X_t\right)^2}\;\1_{\left\{X_t<x\right\}}\left[{\rm
     d}X_t-S\left(\vartheta,X_t\right)\,{\rm d}t\right],\quad x\in {\cal R},
$$
and (formally) the statistic
\begin{align*}
\hat U_T\left(x\right)=U_T\left(x,\hat\vartheta_T ,X^T\right).
\end{align*}
Note that with $\Pb_{\vartheta _0}$ probability 1 we have the equality
$
U_T\left(\infty  ,\vartheta_0,X^T\right)=U_T\left(\vartheta_0  ,X^T\right).
$
The asymptotic behaviour of this statistic can be explained as follows
(again, formally).
\begin{align}
\hat U_T\left(x\right)&=\varphi _T(\hat\vartheta_T
 )\int_{0}^{T}\frac{\dot S(\hat\vartheta_T ,X_t)}{\sigma
  \left(X_t\right)^2}\;\1_{\left\{X_t<x\right\}}\left[{\rm
     d}X_t-S(\hat\vartheta_T,X_t)\,{\rm d}t\right]\nonumber\\
&=\varphi _T(\hat\vartheta_T
 )\int_{0}^{T}\frac{\dot S(\hat\vartheta_T ,X_t)}{\sigma
  \left(X_t\right)}\;\1_{\left\{X_t<x\right\}}{\rm
     d}W_t\nonumber\\
&\quad +\varphi _T(\hat\vartheta_T
 )\int_{0}^{T}\frac{\dot S(\hat\vartheta_T
  ,X_t)\left[S(\vartheta_0,X_t)-S(\hat\vartheta_T,X_t)\right]}{\sigma
  \left(X_t\right)^2}\;\1_{\left\{X_t<x\right\}}\,{\rm d}t\nonumber\\
&=\varphi _T(\hat\vartheta_T
 )\int_{0}^{T}\frac{\dot S(\hat\vartheta_T ,X_t)}{\sigma
  \left(X_t\right)}\;\1_{\left\{X_t<x\right\}}{\rm
     d}W_t\nonumber\\
&\quad -\frac{\hat\vartheta_T-\vartheta _0 }{\varphi _T(\hat\vartheta_T
 )}\int_{0}^{T}\frac{\dot S(\hat\vartheta_T ,X_t)\dot
  S(\tilde\vartheta_T,X_t)}{T {\rm I}\left(\hat\vartheta_T
\right)\sigma
  \left(X_t\right)^2}\;\1_{\left\{X_t<x\right\}}\,{\rm d}t\nonumber\\
&=\varphi _T(\vartheta_0
 )\int_{0}^{T}\frac{\dot S(\vartheta_0 ,X_t)}{\sigma
  \left(X_t\right)}\;\1_{\left\{X_t<x\right\}}{\rm
     d}W_t\nonumber\\
&\quad -\int_{0}^{T}\frac{\dot S(\vartheta_0 ,X_t)}{\sqrt{T{\rm I}\left(\vartheta _0
\right)}\sigma
  \left(X_t\right)}\;{\rm
     d}W_t\int_{0}^{T}\frac{\dot S(\vartheta _0 ,X_t)^2}{T {\rm I}\left(\vartheta _0
\right)\sigma
  \left(X_t\right)^2}\;\1_{\left\{X_t<x\right\}}\,{\rm d}t+o\left(1\right).
\label{eg}
\end{align}
Here $\vartheta _0$ is the true value of the parameter. These integrals have the following limits
\begin{align*}
&\frac{1}{\sqrt{T {\rm I}\left(\vartheta _0
\right)}}\int_{0}^{T}\frac{\dot S(\vartheta_0 ,X_t)}{\sigma
  \left(X_t\right)}\;\1_{\left\{X_t<x\right\}}{\rm
     d}W_t\Longrightarrow \int_{-\infty }^{x}\frac{\dot S\left(\vartheta
  _0,y\right)\sqrt{f\left(\vartheta _0,y\right)}}{\sqrt{{\rm I}\left(\vartheta _0
\right)}\;\sigma \left(y\right)}\;{\rm
  d}w\left(y\right) ,\\
&\frac{1}{T {\rm I}\left(\vartheta _0
\right)}\int_{0}^{T}\frac{\dot S(\vartheta _0 ,X_t)^2}{\sigma
  \left(X_t\right)^2}\;\1_{\left\{X_t<x\right\}}\,{\rm d}t\longrightarrow
  \int_{-\infty }^{x}\frac{\dot S(\vartheta _0 ,y)^2\,{f\left(\vartheta
      _0,y\right)}}{ {\rm I}\left(\vartheta _0
\right)\,\sigma
  \left(y\right)^2}\;{\rm d}y.
\end{align*}
Let us denote
$$
\tau =\int_{-\infty }^{x}\frac{\dot S(\vartheta _0 ,y)^2\,{f\left(\vartheta
      _0,y\right)}}{ {\rm I}\left(\vartheta _0
\right)\,\sigma
  \left(y\right)^2}\;{\rm d}y,\qquad 0\leq \tau \leq 1.
$$
Then we have the convergence
$$
\hat U_T\left(x\right)\Longrightarrow W\left(\tau
\right)-W\left(1\right)\,\tau =B\left(\tau \right),\qquad 0\leq \tau \leq 1.
$$
This limit suggests the construction of the statistic
\begin{align*}
\Delta _T=\int_{-\infty }^{\infty }\frac{\hat U_T\left(x\right)^2\,\dot
  S(\hat\vartheta _T ,x)^2}{{\rm I}(\hat\vartheta _T
)\,\sigma
  \left(x\right)^2} \,{\rm d}F(\hat\vartheta _T,x)
\end{align*}
and the test
$$
\hat\psi _T=\1_{\left\{\Delta _T>c_\alpha \right\}},\qquad \Pb\left(\Delta
>c_\alpha \right)=\alpha.
$$
Note that
$$
\tau _T=\int_{-\infty }^{x}\frac{\dot S(\hat\vartheta _T ,x)^2}{{\rm
    I}(\hat\vartheta _T
)\,\sigma
  \left(x\right)^2} \,{\rm d}F(\hat\vartheta _T,x)\longrightarrow \tau .
$$
Hence if we verify that $\Delta _T\Rightarrow \Delta $, then the test
$\hat\psi _T\in {\cal K}_\alpha $ and is ADF.

We have the same problem with the definition of the stochastic integral
\begin{align*}
\int_{0}^{T}\frac{\dot S(\hat\vartheta_t ,X_t)}{\sigma
  \left(X_t\right)^2}\;\1_{\left\{X_t<x\right\}}\;{\rm
     d}X_t
\end{align*}
as in \eqref{i} and we propose two approaches. In the first one
we replace it  by the ordinary integral  using the It\^o formula as it was done
above and in the second approach we propose using a preliminary consistent
estimator of the parameter $\vartheta $.

 {\it First approach.} Introduce the function,
\begin{align*}
H_T\left(\vartheta ,x,z\right)=\int_{X_0}^{z}\frac{\dot S(\vartheta
  ,y)}{\sigma
  \left(y\right)^2}\phi_T\left(x-y\right){\rm      d}y,\; H\left(\vartheta
,x,z\right)=\int_{X_0}^{z}\frac{\dot S(\vartheta ,y)}{\sigma
  \left(y\right)^2}\1_{\left\{y<x\right\}}{\rm      d}y
\end{align*}
where $\phi_T\left(x-y\right) $ is a ``smooth approximation''  of the
indicator function  $ \1_{\left\{y<x\right\}}$. For example,
$\phi_T\left(x-y\right)=\phi\left(\frac{x-y}{d_T}\right) $, where
\begin{align*}
\phi \left(z\right)=a^{-1}\int_{-\infty
}^{z}e^{\frac{v^2}{v^2-1}}\;\1_{\left\{\left|v\right|<1\right\}}\;{\rm d}v,\quad  a=\int_{-1
}^{1}e^{\frac{v^2}{v^2-1}}\1_{\left\{\left|v\right|<1\right\}}\;{\rm d}v
\end{align*}
and $d_T\rightarrow 0$.

We write
\begin{align*}
\int_{0}^{T}\frac{\dot S(\vartheta ,X_t)}{\sigma
  \left(X_t\right)^2}\;\phi_T\left(x-X_t\right) \;{\rm
  d}X_t&=H_T\left(\vartheta
,x,X_T\right)\\
&\quad -\frac{1}{2}\int_{0}^{T}{\sigma
  \left(X_s\right)^2}{
  }\left(H_T\right)''_{z,z}\left(\vartheta ,x,X_s\right){\rm d}s.
\end{align*}
Then we use the representation of the modified score-function process $\tilde
U_T\left(x\right)$ (we replaced the indicator function by its smooth
approximation)
\begin{align*}
\tilde U_T\left(x\right)&=\varphi _T(\hat\vartheta _T)H_T(\hat\vartheta _T
,x,X_T) -\frac{\varphi _T(\hat\vartheta _T)}{2}\int_{0}^{T}{\sigma
  \left(X_s\right)^2}{
  }\left(H_T\right)''_{z,z}(\hat\vartheta _T ,x,X_s){\rm d}s\\
&\quad -\varphi _T(\hat\vartheta_T
 )\int_{0}^{T}\frac{\dot S(\hat\vartheta_T ,X_t)S(\vartheta_0,X_t)}{\sigma
  \left(X_t\right)^2}\;\phi_T\left(x-X_s\right)\,{\rm d}t\\
&\quad
+\varphi _T(\hat\vartheta_T
 )\int_{0}^{T}\frac{\dot S(\hat\vartheta_T
  ,X_t)\left[S(\vartheta_0,X_t)-S(\hat\vartheta_T,X_t)\right]}{\sigma
  \left(X_t\right)^2}\;\phi_T\left(x-X_s\right)\,{\rm d}t\\
&=J_T\left(\hat\vartheta _T,x\right)-K_T\left(\hat\vartheta _T,x\right).
\end{align*}

Direct but cumbersome calculations give   the limits
\begin{align*}
&J_T\left(\hat\vartheta _T,x\right)\Longrightarrow \int_{-\infty }^{x
  }\frac{\dot S\left(\vartheta
  _0,y\right)\sqrt{f\left(\vartheta _0,y\right)}}{\sqrt{{\rm I}\left(\vartheta _0
\right)}\;\sigma \left(y\right)}\;{\rm
  d}w\left(y\right),\\
&K_T\left(\hat\vartheta _T,x\right)= \frac{\left(\hat\vartheta_T-\vartheta
    \right)}{\varphi _T(\hat\vartheta_T
 )  }
\int_{0}^{T}\frac{\dot S(\hat\vartheta_T
  ,X_t)\dot S(\tilde\vartheta_T,X_t)}{T {\rm I}\left(\hat\vartheta _T\right)\sigma
  \left(X_t\right)^2}\;\1_{\left\{X_t<x\right\}})\,{\rm
  d}t\left(1+o\left(1\right)\right)\\
&\qquad \quad \Longrightarrow  \int_{-\infty }^{\infty }\frac{\dot S\left(\vartheta
  _0,y\right)\sqrt{f\left(\vartheta _0,y\right)}}{\sqrt{{\rm I}\left(\vartheta _0
\right)}\;\sigma \left(y\right)}\;{\rm
  d}w\left(y\right)\;\int_{-\infty }^{x}\frac{\dot S\left(\vartheta
  _0,y\right)\sqrt{f\left(\vartheta _0,y\right)}}{\sqrt{{\rm I}\left(\vartheta _0
\right)}\;\sigma \left(y\right)}\;{\rm
  d}y.
\end{align*}
Thus we have the following result.
\begin{proposition}
\label{P3} Suppose that the conditions of regularity are fulfilled, then the
test $\tilde\psi _T =\1_{\left\{\tilde\Delta _T >c_\alpha
  \right\}}$ with
\begin{align*}
\tilde\Delta _T=\int_{-\infty  }^{\infty }\frac{\tilde U_T
  \left(x\right)^2\dot S(\hat\vartheta _{T
  },x)^2}{{\rm I}\left(\hat\vartheta _{T
  }  \right)\sigma \left(x\right)^2}\,{\rm d}F(\hat\vartheta _{T
  },x )
\end{align*}
is ADF and belongs to ${\cal K}_\alpha $.
\end{proposition}
{\it The second approach.} Let us introduce a consistent preliminary estimator
$\bar\vartheta _{\sqrt{T}}$ constructed using the first $X_t,0\leq t\leq
\sqrt{T}$ observations. For example,  the method of moments estimator can be used (see conditions of consistency in \cite{Kut04}, Section 2.4). The
corresponding statistic  is
\begin{align*}
 V_T\left(x\right)&=\varphi _T(\bar\vartheta_{\sqrt{T}}
 )\int_{\sqrt{T}}^{T}\frac{\dot S(\bar\vartheta_{\sqrt{T}} ,X_t)}{\sigma
  \left(X_t\right)^2}\;\1_{\left\{X_t<x\right\}}\left[{\rm
     d}X_t-S(\hat\vartheta_T,X_t)\,{\rm d}t\right].
\end{align*}
The stochastic integral is well defined and its limit can be obtained
 calculations, similar to \eqref{eg}.

\begin{proposition}
\label{P4} Suppose that the conditions of regularity are fulfilled and the
preliminary estimator $\bar\vartheta _{\sqrt{T}}  $ is consistent, then the
test $\tilde\psi _T =\1_{\left\{\tilde\Delta _T >c_\alpha
  \right\}}$ with
\begin{align*}
\tilde\Delta _T=\int_{-\infty  }^{\infty }\frac{V_T
  \left(x\right)^2\dot S(\bar\vartheta _{T
  },x)^2}{{\rm I}\left(\bar\vartheta _{T
  }  \right)\sigma \left(x\right)^2}\,{\rm d}F(\bar\vartheta _{T
  },x )
\end{align*}
is ADF and belongs to ${\cal K}_\alpha $.
\end{proposition}

The condition of the consistency is
\begin{align*}
\inf_{\bar\vartheta ,\hat\vartheta }\int_{-\infty }^{\infty }\frac{M(\bar\vartheta ,\hat\vartheta
,x)^2\;\dot S\left(\bar\vartheta
  ,x\right)^2}{\sigma \left(x\right)^2}\;{f\left(x\right)}\;{\rm d}x>0,
\end{align*}
where
$$
M(\bar\vartheta ,\hat\vartheta ,x)=\int_{-\infty }^{x}\frac{\dot S\left(\bar\vartheta
  ,y\right) \left[S\left(y\right)-S{\left(\hat\vartheta
      ,y\right)}\right]{}}{\sigma \left(y\right)^2}\;f\left(y\right)\;{\rm
  d}y.
$$
\subsection{Periodic Poisson Processes.}
The last observations model  is  a periodic Poisson process
$$
X^n=\left(X_t,0\leq t\leq
T=n\tau_*\right)
$$
of known period $\tau_* >0$.
For $0\leq s<t$ and $k=0,1,2,\ldots $
$$
\Pb\left(X_t-X_s=k\right)=\frac{\left[\Lambda \left(t\right)-\Lambda
    \left(s\right)\right]^k}{k!}  \exp \left\{-\Lambda \left(t\right)+\Lambda
    \left(s\right)\right\}.
$$
The mean $\Lambda \left(t\right)$ and
intensity function $\lambda \left(t\right)$ satisfy the relations
$$
\Lambda \left(t\right)=\Ex X_t, \qquad \Lambda
\left(t\right)=\int_{0}^{t}\lambda \left(s\right)\,{\rm d}s
$$
and $\lambda \left(t+k\tau_* \right)=\lambda \left(t \right)$.

We observe a trajectory $X^n$ of the Poisson process of intensity  function
$\lambda \left(\cdot \right)$ and we  have to test the hypothesis

${\scr H}_0$ : {\it The intensity function  $\lambda \left(t\right)=\lambda
  \left(\vartheta
,t\right),\, \vartheta \in \Theta=\left(a,b\right)$. }

{Here $\lambda \left(\vartheta ,\cdot \right)$ is some known   function
  satisfying the following conditions of regularity}.

 The intensity function $\lambda \left(\vartheta ,\cdot \right)$ is twice
 continuously differentiable w.r.t. $\vartheta $, strictly positive and the
 identifiability condition   holds: for any $\nu >0$
$$
\inf_{\vartheta _0\in \Theta }\inf_{\left|\vartheta
  -\vartheta _0\right|>\nu } \int_{0}^{\tau_*
}\left[\sqrt{\lambda \left(\vartheta ,s\right)}-\sqrt{\lambda
    \left(\vartheta_0 ,s\right)}\right]^2 \,{\rm d}s>0 .
$$
 The likelihood ratio function is
\begin{align*}
L\left(\vartheta ,X^n\right)=\exp\left\{\sum_{j=1}^{n}\int_{0}^{\tau _*}\ln
{\lambda\left(\vartheta ,t\right) }{}\;{\rm d}X_j\left(t\right)-n\int_{0}^{\tau
  _*}\left[\lambda \left(\vartheta ,t\right)-1\right]\; {\rm d}t\right\}
\end{align*}
and the MLE $\hat\vartheta _n$ is defined by the equation like \eqref{1.3}.
Then the MLE  admits the representation
$$
\sqrt{n}\left(\hat\vartheta _T-\vartheta \right)=\frac{1}{{\rm I}\left(\vartheta
\right)\sqrt{n}}\sum_{j=1}^{n} \int_{0}^{\tau_*  }\frac{\dot \lambda \left(\vartheta
  ,s\right)}{\lambda \left(\vartheta
  ,s\right)} \; {\rm d}\left[X_j\left(s\right)-\lambda \left(\vartheta
  ,s\right){\rm d}s\right]   +o\left(1\right).
$$
Here  $X_j\left(s\right)=X_{\left(j-1\right)\tau_*
    +s}-X_{\left(j-1\right)\tau_*  },\:0\leq s\leq \tau_* ,\:  j=1,2,\ldots,n$ and
$ {\rm I}\left(\vartheta
\right)$  is the Fisher information
$$
{\rm I}\left(\vartheta
\right)=\int_{0 }^{\tau_*
}\frac{\dot \lambda \left(\vartheta ,s\right)^2}{\lambda
    \left(\vartheta ,s\right)}\,{\rm d}s>0.
$$
The proof can be found in \cite{Kut98}.

The score-function for this process is
\begin{align*}
\frac{\partial \ln L\left(\vartheta ,X^n\right)}{\partial \vartheta
}=\sum_{j=1}^{n}\int_{0}^{\tau _*}\frac {\dot \lambda\left(\vartheta ,t\right)
}{\lambda\left(\vartheta ,t\right)}\;\left[{\rm
    d}X_j\left(t\right)-\lambda\left(\vartheta ,t\right)\,{\rm d}t\right]
\end{align*}
and we define the normalized score-function process
\begin{align*}
U_n\left(t,\vartheta ,X^n\right)=\frac{1}{\sqrt{{\rm I}\left(\vartheta
  \right)n}}\sum_{j=1}^{n}\int_{0}^{t}\frac {\dot
  \lambda\left(\vartheta ,s\right) }{\lambda\left(\vartheta
  ,s\right)}\;\left[{\rm d}X_j\left(s\right)-\lambda\left(\vartheta
  ,s\right)\,{\rm d}s\right].
\end{align*}
We construct  the GoF test with  the help of the statistic
$$
\hat U_n\left(t\right)=U_n\left(t,\hat\vartheta _n,X^n\right)
$$
Its formal expansion  provides us with the following expressions (we put below $\pi
_j\left(s\right)=X_j\left(s\right)-\Lambda \left(\vartheta ,s\right)$)
 \begin{align*}
\hat U_n\left(t\right)&=\frac{1}{\sqrt{{\rm I}(\hat\vartheta _n
  )n}}\sum_{j=1}^{n}\int_{0}^{t}\frac {\dot
  \lambda (\hat\vartheta _n ,s ) }{\lambda (\hat\vartheta _n
  ,s )}\;\left[{\rm d}X_j\left(s\right)-\lambda (\hat\vartheta _n
  ,s )\,{\rm d}s\right]\\
&=\frac{1}{\sqrt{{\rm I}(\hat\vartheta _n
  )n}}\sum_{j=1}^{n}\int_{0}^{t}\frac {\dot
  \lambda (\hat\vartheta _n ,s ) }{\lambda (\hat\vartheta _n
  ,s )}\;\left[{\rm d}X_j\left(s\right)-\lambda (\vartheta
  ,s )\,{\rm d}s\right]\\
&\qquad +\frac{1}{\sqrt{{\rm I}(\hat\vartheta _n
  )n}}\sum_{j=1}^{n}\int_{0}^{t}\frac {\dot
  \lambda (\hat\vartheta _n ,s ) }{\lambda (\hat\vartheta _n
  ,s )}\;\left[\lambda (\vartheta
  ,s )-\lambda (\hat\vartheta _n
  ,s )\right]\,{\rm d}s\\
&=\frac{1}{\sqrt{{\rm I}(\vartheta
  )n}}\sum_{j=1}^{n}\int_{0}^{t}\frac {\dot
  \lambda (\vartheta ,s ) }{\lambda (\vartheta
  ,s )}\;{\rm d}\pi _j\left(s\right)-\frac{\sqrt{n}(\hat\vartheta
  _n-\vartheta)}{\sqrt{{\rm I}(\vartheta
  )}}\int_{0}^{t}\frac {\dot
  \lambda (\vartheta,s )^2 }{\lambda (\vartheta
  ,s )}\,{\rm d}s+o\left(1\right)\\
&=\frac{1}{\sqrt{{\rm I}(\vartheta
  )n}}\sum_{j=1}^{n}\int_{0}^{t}\frac {\dot
  \lambda (\vartheta ,s ) }{\lambda (\vartheta
  ,s )}\;{\rm d}\pi _j\left(s\right)\\
&\qquad  -\frac{1}{\sqrt{{\rm I}(\vartheta
  )n}}\sum_{j=1}^{n}\int_{0}^{\tau _*}\frac {\dot
  \lambda (\vartheta ,s ) }{\lambda (\vartheta
  ,s )}\;{\rm d}\pi _j\left(s\right)\;\int_{0}^{t}\frac {\dot
  \lambda (\vartheta,s )^2 }{{\rm I}(\vartheta
  )\,\lambda (\vartheta
  ,s )}\,{\rm d}s+o\left(1\right).
\end{align*}
By the central limit theorem we have the convergence in distribution
\begin{align*}
&\frac{1}{\sqrt{{\rm I}(\vartheta
  )n}}\sum_{j=1}^{n}\int_{0}^{t}\frac {\dot
  \lambda (\vartheta ,s ) }{\lambda (\vartheta
  ,s )}\;{\rm d}\pi _j\left(s\right)\Longrightarrow \frac{1}{\sqrt{{\rm I}(\vartheta
  )}}\int_{0}^{t}\frac {\dot
  \lambda (\vartheta ,s ) }{\sqrt{\lambda (\vartheta
  ,s )}}\;{\rm d}W_s,\\
&\frac{1}{\sqrt{{\rm I}(\vartheta
  )n}}\sum_{j=1}^{n}\int_{0}^{\tau _*}\frac {\dot
  \lambda (\vartheta ,s ) }{\lambda (\vartheta
  ,s )}\;{\rm d}\pi _j\left(s\right)\Longrightarrow \frac{1}{\sqrt{{\rm I}(\vartheta
  )}}\int_{0}^{\tau _*}\frac {\dot
  \lambda (\vartheta ,s ) }{\sqrt{\lambda (\vartheta
  ,s )}}\;{\rm d}W_s,
\end{align*}
where $W_t,0\leq t\leq \tau _*$ is some Wiener process.
Therefore, if we put
$$
 \tau =\int_{0}^{t}\frac {\dot \lambda (\vartheta,s )^2 }{{\rm I}(\vartheta
  )\,\lambda (\vartheta ,s )}\,{\rm d}s,\qquad 0\leq \tau \leq 1,
$$
then once again we obtain the convergence
$$
\hat U_n\left(t\right)\Longrightarrow W\left(\tau \right)-W\left(1 \right)\tau
=B\left(\tau \right),\qquad 0\leq \tau \leq 1.
$$
We can consider two approaches as before, but we present here the
second (more simple)   construction of the test. Let us take any
consistent estimator $\bar\vartheta _{{N}}$ of the parameter $\vartheta $
constructed by the first $N=\left[\sqrt{n}\right]$ observations
$X^N=\left(X_1,\ldots,X_N\right)$.  Then we set
\begin{align*}
V_n\left(t\right)&=\frac{1}{\sqrt{{\rm I}(\bar\vartheta _N
  )n}}\sum_{j=N+1}^{n}\int_{0}^{t}\frac {\dot
  \lambda (\bar\vartheta _N ,s ) }{\lambda (\bar\vartheta _N
  ,s )}\;\left[{\rm d}X_j\left(s\right)-\lambda (\hat\vartheta _n
  ,s )\,{\rm d}s\right].
\end{align*}
The estimator $\bar\vartheta _N$ and the observations
$X_{N+1}^n=\left(X_{N+1},\ldots,X_n\right)$ are independent and the stochastic
integral with respect to the Poisson process is well defined (see Liptser,
Shiryayev \cite{LS}, Section 18.4).

\begin{proposition}
\label{P5}
Let  the conditions of regularity  be fulfilled, then the test $\tilde \psi
_n=\1_{\left\{\Delta _n>c_\alpha \right\}}$ with
$$
\Delta _n=\int_{0}^{\tau _*}\frac {V_n\left(t\right)^2\dot \lambda
  (\bar\vartheta _N,s )^2 }{{\rm I}(\bar\vartheta _N )\,\lambda (\bar\vartheta _N ,s
  )}\,{\rm d}s.
$$
is ADF  and  belongs to ${\cal K}_\alpha $.
\end{proposition}
To prove this proposition we have to verify the convergence
$$
\Delta _n\Longrightarrow \Delta =\int_{0}^{1}B\left(\tau \right)^2{\rm d}\tau
$$
under hypothesis ${\cal   H}_0 $

{\bf Example.} Suppose that the intensity function under hypothesis ${\cal
  H}_0 $ is
\begin{align*}
\lambda \left(\vartheta ,t\right)=\vartheta h\left(t\right)+\lambda _0,\quad 0\leq t\leq \tau _*,
\end{align*}
where $\vartheta \in \Theta =\left(a,b\right), a>0$ and the function $h\left(t\right)>0$.

Then we can take as preliminary estimator the minimum distance estimator
\begin{align*}
\bar\vartheta _N&=\arg\inf_{\vartheta\in \Theta } \int_{0}^{\tau
  _*}\left[\hat\Lambda _N\left(t\right) -\vartheta H\left(t\right)-\lambda
  _0t\right]^2{\rm d}t \\
&=\frac{ \int_{0}^{\tau
  _*}\left[\hat\Lambda _N\left(t\right) -\lambda
  _0t\right]H\left(t\right)  {\rm d}t}{\int_{0}^{\tau
  _*}H\left(t\right)^2{\rm d}t}.
\end{align*}
Here
\begin{align*}
\hat\Lambda _N\left(t\right)=\frac{1}{N}\sum_{j=1}^{N}X_j\left(t\right),\qquad
H\left(t\right)=\int_{0}^{t}h\left(s\right)\,{\rm d}s .
\end{align*}
This is an unbiased, consistent and asymptotically normal estimator of the
parameter $\vartheta $.

The score-function process $V_n\left(\cdot \right) $ and the test statistics $\tilde \Delta _n
$ are
\begin{align*}
V_n\left(t\right)&=\frac{1}{\sqrt{{\rm I}(\bar\vartheta_{{N}}
  )n}}\sum_{j=N+1}^{n}\int_{0}^{t}\frac {h\left(s\right) }{\bar\vartheta_{{N}}
  h\left(s\right)+\lambda _0}\;\left[{\rm d}X_j\left(s\right)-\left[\bar\vartheta_{{N}}
  h\left(s\right)+\lambda
  _0\right]\,{\rm d}s\right],\\
\tilde \Delta _n&=\int_{0}^{\tau _*}\frac {V_n\left(t\right)^2h\left(s\right)^2   }{{\rm I}(\bar\vartheta _N )\,\left[\bar\vartheta_{{N}}
  h\left(s\right)+\lambda
  _0\right]}\,{\rm d}s,\qquad {\rm I}\left(\vartheta \right)=\int_{0}^{\tau
  _*}\frac{h\left(t\right)^2}{\vartheta h\left(t\right)+\lambda _0}\;{\rm d}t ,
\end{align*}
respectively.

\section{Other tests and models}

\subsection{Other tests}

The statistics $U_\varepsilon \left(\cdot \right),U_T \left(\cdot \right)$ and
$U_n \left(\cdot \right) $ can be used for  construction of the
ADF GoF tests of Kolmogorov-Smirnov type. For example, the following convergence
\begin{align*}
\Delta _\varepsilon ^*=\sup_{\nu _\varepsilon \leq t\leq T}\left|V_\varepsilon
\left(t\right)\right|\Longrightarrow
\sup_{0\leq \tau \leq 1}\left|B\left(\tau \right)\right|=\Delta ^*
\end{align*}
can be easily proved.
Hence the test
\begin{align*}
\psi _\varepsilon ^*=\1_{\left\{\Delta _\varepsilon ^*> d_\alpha
  \right\}},\qquad \Pb\left(\Delta ^* >d_\alpha  \right)=\alpha
\end{align*}
belongs to ${\cal K}_\alpha $ and is ADF.
Of course  similar tests can be constructed in the cases of observations of
the  ergodic
diffusion and inhomogeneous Poisson processes as well.

\subsection{Nonlinear AR process}

Suppose that the observations $X^n=\left(X_0,X_1,\ldots,X_n\right)$ satisfy
the relation
\begin{align*}
X_j=S\left(X_{j-1}\right)+\varepsilon _j,\quad j=1,\ldots,n
\end{align*}
 and we have to test a parametric hypothesis
\begin{align*}
{\cal H}_0\qquad :\qquad S\left(x\right)=S\left(\vartheta ,x\right),\quad
\vartheta \in \Theta =\left(a,b\right).
\end{align*}
Here $S\left(\vartheta ,x\right)$ is some known function and $\vartheta $ is the
unknown parameter. The random variables $\varepsilon _1,\ldots,\varepsilon _j$
are i.i.d. with the known density function $f\left(x\right)$.

The functions $S\left(\vartheta ,x\right)$ and $f\left(x\right)>0$ are such that
the time series  $\left(X_j\right)_{j\geq 1}$ has ergodic properties with the
density of invariant law $\varphi \left(\vartheta ,x\right)$ for all
$\vartheta \in \Theta $, i.e., for any
function $h\left(\cdot \right)$ such that  $\Ex_\vartheta \left|h\left(\xi
\right)\right|<\infty $ (here $\xi \sim \varphi \left(\vartheta ,\cdot
\right)$) we have the law of large
numbers
\begin{align*}
\frac{1}{n}\sum_{j=1}^{n}h\left(X_j\right)\longrightarrow \Ex_\vartheta
h\left(\xi \right).
\end{align*}
Moreover we suppose that the tails of $\varphi \left(\vartheta ,x\right)$
decrease sufficiently fast
\begin{equation}
\label{tails}
\varphi \left(\vartheta ,x\right)\leq \frac{C}{\left|x\right|^{1+\gamma} }
\end{equation}
with some positive constants $\gamma$ and $ C$, which do not depend on $\vartheta $.
The log-density function $\ell\left(x\right)=\ln f\left(x\right)$ has three
continuous bounded derivatives $\ell'\left(x \right),$ $\ell''\left(x \right),$ $\ell'''\left(x
\right)$ and the function $S\left(\vartheta ,x\right)$ has two continuous
bounded derivatives $\dot S\left(\vartheta ,x\right),$ $\ddot S\left(\vartheta
,x\right)$  w.r.t. $\vartheta $.

The log-likelihood  function is
\begin{align*}
L\left(\vartheta ,X^n\right)=\ln\varphi \left(\vartheta
,X_0\right)+\sum_{j=1}^{n}\ln f\left(X_j-S\left(\vartheta
,X_{j-1}\right)\right),\qquad \vartheta \in \left(a,b\right).
\end{align*}
We suppose that the initial value $X_0$ has invariant density function
$\varphi \left(\vartheta ,x\right)$ and therefore the time series $\left(X_j
\right)_{j\geq 0}$ is stationary.

The Score-function is
\begin{align*}
U_n\left(\vartheta ,X^n\right)=-\sum_{j=1}^{n} \ell'\left(X_j-S\left(\vartheta
,X_{j-1}\right) \right)\dot S\left(\vartheta ,X_{j-1}\right).
\end{align*}

Also we assume  that the regularity conditions are fulfilled so that the MLE
$\hat\vartheta _n$ is consistent and admits the representation
\begin{align}
\label{mle}
\sqrt{n}\left(\hat\vartheta _n-\vartheta \right)=\frac{-1}{{\rm
    I}\left(\vartheta \right)\sqrt{n}} \sum_{j=1}^{n} \ell'\left(X_j-S\left(\vartheta ,X_{j-1}\right) \right)\dot S\left(\vartheta ,X_{j-1}\right)+o\left(1\right),
\end{align}
where the Fisher information
\begin{align*}
{\rm I}\left(\vartheta \right)=\EE_\vartheta \left[\ell'\left(\varepsilon _1
  \right)\dot S\left(\vartheta ,\xi \right)\right]^2={\rm I}_f\;{\rm
  I}_\vartheta.
\end{align*}
Here we denoted $\EE_\vartheta $ the expectation related to the couple of
independent random variables $\left(\varepsilon ,\xi \right)$, i.e.,
\begin{align*}
 {\rm I}_f =\Ex \ell'\left(\varepsilon
\right)^2=\int_{-\infty }^{\infty
}\frac{f'\left(x\right)^2}{f\left(x\right)}{\rm d}x,\quad {\rm
  I}_\vartheta=\Ex_\vartheta \dot S\left(\vartheta ,\xi \right)^2 =\int_{-\infty }^{\infty
}{\dot S\left(\vartheta ,x\right)^2}{\varphi \left(\vartheta ,x\right)}{\rm d}x.
\end{align*}
Note that from this representation and the central limit theorem it follows
that the MLE is asymptotically normal (see, e.g., \cite{HB93})
\begin{align*}
\hat u_n=\sqrt{n}\left(\hat\vartheta _n-\vartheta \right)\Longrightarrow {\cal N}\left(0,{\rm
    I}\left(\vartheta \right)^{-1}\right).
\end{align*}
Introduce the normalized score-function process
\begin{align*}
U_n\left(x,\vartheta ,X^n\right)=\frac{-1}{\sqrt{{\rm
    I}\left(\vartheta \right)n}  }\sum_{j=1}^{n} \ell'\left(X_j-S\left(\vartheta
,X_{j-1}\right) \right)\dot S\left(\vartheta ,X_{j-1}\right)\,\1_{\left\{X_{j-1<x}\right\}}
\end{align*}
and the corresponding statistics
\begin{align*}
\hat U_n\left(x\right)=\frac{-1}{\sqrt{{\rm
    I}(\hat\vartheta_n )n}  }\sum_{j=1}^{n} \ell'\left(X_j-S(\hat\vartheta_n
,X_{j-1}) \right)\dot S(\hat\vartheta_n ,X_{j-1})\,\1_{\left\{X_{j-1<x}\right\}}.
\end{align*}
Using the expansion at the vicinity of the true value $\vartheta $ we can write
\begin{align*}
&\hat U_n\left(x\right)=\frac{-1}{\sqrt{{\rm
    I}(\vartheta )n}  }\sum_{j=1}^{n} \left[\ell'\left(X_j-S\left(\vartheta
,X_{j-1}\right) \right) \right.\\
&\qquad \left.-\frac{\hat u_n}{\sqrt{n}}\ell''\left(X_j-S\left(\vartheta
,X_{j-1}\right) \right)\dot S\left(\vartheta ,X_{j-1}\right) \right]\dot
  S\left(\vartheta
  ,X_{j-1}\right)\,\1_{\left\{X_{j-1<x}\right\}}+o\left(1\right)\\
&\quad =\frac{-1}{\sqrt{{\rm
    I}(\vartheta )n}  }\sum_{j=1}^{n} \ell'\left(X_j-S\left(\vartheta
,X_{j-1}\right) \right)\dot S\left(\vartheta ,X_{j-1}\right)\,\1_{\left\{X_{j-1<x}\right\}}\\
&\qquad+
\frac{{\hat u_n}\sqrt{{\rm     I}(\vartheta )}}{{\rm
    I}(\vartheta )n}
\sum_{j=1}^{n} \ell''\left(X_j-S\left(\vartheta
,X_{j-1}\right) \right)
\dot S\left(\vartheta ,X_{j-1}\right)^2\,\1_{\left\{X_{j-1<x}\right\}}+o\left(1\right).
\end{align*}
The standard arguments allow us to write
\begin{align*}
&\frac{-1}{{\rm    I}(\vartheta )n}
\sum_{j=1}^{n} \ell''\left(X_j-S\left(\vartheta
,X_{j-1}\right) \right)
\dot S\left(\vartheta ,X_{j-1}\right)^2\,\1_{\left\{X_{j-1<x}\right\}}\\
& \quad =\frac{-1}{{\rm    I}(\vartheta )n}
\sum_{j=1}^{n} \ell''\left(\varepsilon _j\right)
\dot S\left(\vartheta ,X_{j-1}\right)^2\,\1_{\left\{X_{j-1<x}\right\}}
\longrightarrow\frac{1}{{\rm    I}_\vartheta   } \int_{-\infty }^{x
}{\dot S\left(\vartheta ,y\right)^2}{\varphi \left(\vartheta ,y\right)}{\rm d}y.
\end{align*}
Recall, that
\begin{align*}
\Ex \ell''\left(\varepsilon  \right)=\Ex\left(\frac{f''\left(\varepsilon
  \right)f\left(\varepsilon \right)-f'\left(\varepsilon
  \right)^2}{f\left(\varepsilon \right)^2}\right) =-\Ex\left(\frac{f'\left(\varepsilon
  \right)}{f\left(\varepsilon \right)}\right)^2=-{\rm I}_f
\end{align*}
because
$$
\Ex\left(\frac{f''\left(\varepsilon
  \right)}{f\left(\varepsilon \right)}\right)=\int_{-\infty }^{\infty }f''\left(y  \right)\,{\rm d}y=0.
$$
Let us denote
\begin{align*}
W_n\left(x\right) =\frac{-1}{\sqrt{{\rm
    I}(\vartheta )n}  }\sum_{j=1}^{n} \ell'\left(X_j-S\left(\vartheta
,X_{j-1}\right) \right)\dot S\left(\vartheta ,X_{j-1}\right)\,\1_{\left\{X_{j-1<x}\right\}}.
\end{align*}
We have
\begin{align*}
&\Ex_\vartheta W_n\left(x\right)W_n\left(y\right)\\
&=\frac{1}{{{\rm I}(\vartheta
    )n} }\sum_{j=1}^{n} \sum_{i=1}^{n} \Ex_\vartheta \ell'\left(\varepsilon
_j\right)\ell'\left(\varepsilon _i\right) \dot S\left(\vartheta
,X_{j-1}\right)\dot S\left(\vartheta
,X_{i-1}\right)\1_{\left\{X_{j-1<x}\right\}}\,\1_{\left\{X_{i-1<x}\right\}}\\
&=\frac{1}{{{\rm I}(\vartheta
    )} }\EE_\vartheta \ell'\left(\varepsilon
_1\right)^2 \dot S\left(\vartheta
,\xi \right)^2\1_{\left\{\xi <x\wedge y\right\}}\\
&=\min\left({\rm I}_\vartheta ^{-1}\int_{-\infty }^{x}\dot S\left(\vartheta
,z \right)^2\varphi \left(\vartheta ,z\right){\rm d}z,{\rm I}_\vartheta ^{-1}\int_{-\infty }^{y}\dot S\left(\vartheta
,z \right)^2\varphi \left(\vartheta ,z\right){\rm d}z\right)\\
&=\min\left(\tau _x,\tau _y\right),\qquad 0\leq \tau _x={\rm I}_\vartheta ^{-1}\int_{-\infty }^{x}\dot S\left(\vartheta
,z \right)^2\varphi \left(\vartheta ,z\right){\rm d}z\leq 1.
\end{align*}
It can be shown that by the central limit theorem the finite-dimensional
distributions of the random function
$W_n\left(x\right),x\in\RR$ converge to the finite-dimensional distributions
of the Wiener process $W\left(\tau _x\right),x\in \RR$. Moreover the following
estimate holds
\begin{align}
\label{cont}
\Ex_\vartheta\left|W_n\left(x\right)-W_n\left(x\right)\right|^2\leq C\left|x-y\right|.
\end{align}
We have  similar convergence for the MLE due to the representation \eqref{mle}
\begin{align*}
\sqrt{n{\rm     I}(\vartheta )}\left({\hat \vartheta _n-\vartheta }\right)=W_n\left(\infty
\right)+o\left(1\right)\Longrightarrow W\left(1\right)
\end{align*}
with the same Wiener process, i.e., we have the joint asymptotic normality of
$W_n\left(\cdot \right)$ and $\hat u_n$.
Therefore the random functions $\hat U_n\left(x\right)$ have the corresponding
limit
\begin{align*}
\hat U_n\left(x\right)\Longrightarrow W\left(\tau
_x\right)-W\left(1\right)\tau _x=B\left(\tau _x\right)
\end{align*}
and again, we obtain the Brownian bridge $B\left(\tau \right),0\leq \tau \leq 1$.

Let us introduce the statistics
\begin{align*}
\Delta _n=\int_{-\infty }^{\infty }\frac{\hat U_n\left(x\right)^2\dot S(\hat \vartheta _n
,x )^2}{{\rm I}_{\hat \vartheta _n}  }\varphi(\hat \vartheta _n ,x)\,{\rm d}x.
\end{align*}
The convergence of finite-dimensional distributions, the estimate \eqref{cont}
and the condition \eqref{tails} allow us to verify the convergence
\begin{align*}
\Delta _n\Longrightarrow \int_{-\infty }^{\infty }\frac{B\left(\tau _x\right)^2\dot S\left(\vartheta
,x \right)^2}{{\rm I}_\vartheta  }\varphi \left(\vartheta ,x\right){\rm
  d}x=\int_{0}^{1}B\left(\tau \right)^2{\rm d}\tau  .
\end{align*}
Therefore we have the following result.
\begin{proposition}
\label{P6} The test $\hat\psi _n=\1_{\left\{\Delta _n>c_\alpha \right\}}$ is
ADF and belongs to the class ${\cal K}_\alpha $.
\end{proposition}

\bigskip

{\bf Example.} Suppose that the observed time series $\left(X_j\right)_{j\geq
  1}$ under the hypothesis ${\scr H}_0$ is linear AR
\begin{align*}
X_j=\vartheta X_{j-1}+\varepsilon _j,\quad j=1,\ldots,n,
\end{align*}
where  $\vartheta \in \Theta =\left(-1,1\right)$ and $\left(\varepsilon
_j\right)_{j\geq 1} $ are i.i.d. ${\cal N}\left(0,\sigma ^2\right)$
r.v's. Then the aforementioned conditions are satisfied  with the density of
invariant law
 $$
\varphi \left(\vartheta ,x\right)\sim {\cal
  N}\left(0, \frac{\sigma ^2}{1-\vartheta ^2 }\right)
$$
and we assume that $X_0\sim \varphi \left(\vartheta ,x\right)$.

The derivative $\dot S\left(\vartheta ,x\right)=x$ is not bounded, but the
tails of $\varphi \left(\vartheta ,x\right)$ are exponentially decreasing and
the proof of the convergence given above remains valid.

The score-function process is
\begin{align*}
U_n\left(x,\vartheta ,X^n\right)=\frac{1}{\sqrt{n\left(1-\vartheta ^2\right)}}\sum_{j=1}^{n} \left(X_j-\vartheta X_{j-1}\right) X_{j-1}\1_{\left\{X_{j-1}<x\right\}},
\end{align*}
because
$$ {\rm I}\left(\vartheta \right)={\rm I}_f{\rm I}_\vartheta =\frac{1}{\sigma
  ^2}\;\frac{\sigma ^2}{1-\vartheta ^2}=\frac{1}{1-\vartheta ^2}
$$
and we put
\begin{align*}
\hat U_n\left(x\right)=U_n\left(x,\hat\vartheta_n ,X^n\right),\qquad \hat\vartheta _n=\frac{\sum_{j=1}^{n}X_jX_{j-1}}{\sum_{j=1}^{n}X_{j-1}^2}.
\end{align*}
Introduce the statistic
\begin{align*}
\Delta _n= \frac{1-\hat\vartheta _n^2}{\sigma ^2}\int_{-\infty }^{\infty }\hat U_n\left(x\right)^2x^2\varphi(\hat\vartheta _n,x)\,{\rm d}x.
\end{align*}
As it follows from the Proposition \ref{P6}
\begin{align*}
\Delta _n\Longrightarrow \int_{0}^{1}B\left(\tau  \right)^2\,{\rm d}\tau
\end{align*}
and the test $\hat\psi _n=\1_{\left\{\Delta _n>c_\alpha \right\}}$ is ADF and
belongs to ${\cal K}_\alpha $.

\subsection{The case of i.i.d. observations}

Let us see what happens if we apply the
same approach in the case  of i.i.d. observations
$X^n=\left(X_1,\ldots,X_n\right)$,  where $X_j$ has the density function
$f\left(x\right)$. Suppose that we have a parametric hypothesis
\begin{align*}
{\scr H}_0,\qquad :\qquad f\left(x\right)=f\left(\vartheta ,x\right),\quad
\vartheta \in \Theta =\left(a,b\right).
\end{align*}
Here $f\left(\vartheta ,x\right)$ is some known density function satisfying
the regularity conditions, which validate the
calculations below.

 The normalized
score-function statistic is
\begin{align*}
U_n\left(\vartheta ,X^n\right)&=\frac{1}{\sqrt{{\rm I}\left(\vartheta
  \right)n}}\sum_{j=1}^{n}\dot \ell \left(\vartheta
,X_j\right)=\frac{\sqrt{n}}{\sqrt{{\rm I}\left(\vartheta
  \right)}}\int_{-\infty }^{\infty }\dot \ell \left(\vartheta ,y\right)\,{\rm
  d}\hat F_n\left(y\right)\\
&=\frac{\sqrt{n}}{\sqrt{{\rm I}\left(\vartheta
  \right)}}\int_{-\infty }^{\infty }\dot \ell \left(\vartheta ,y\right)\,\left[ {\rm
  d}  \hat F_n\left(y\right)-f\left(\vartheta ,y\right){\rm d}y\right],
\end{align*}
where $\ell \left(\vartheta ,y\right)=\ln f\left(\vartheta ,y\right) $, ${\rm
  I}\left(\vartheta   \right) $ is the Fisher information  and we
used the equality
$$
\int_{-\infty }^{\infty }\dot \ell \left(\vartheta ,y\right)\,f\left(\vartheta
,y\right){\rm d}y=0.
$$
Introduce the score-function process
\begin{align*}
U_n\left(\vartheta ,x,X^n\right)&=\frac{\sqrt{n}}{\sqrt{{\rm I}\left(\vartheta
  \right)}}\int_{-\infty }^{x }\dot \ell \left(\vartheta ,y\right)\,\left[ {\rm
  d}  \hat F_n\left(y\right)-f\left(\vartheta ,y\right){\rm d}y\right],\qquad
x\in {\cal R}
\end{align*}
and the corresponding statistic
\begin{align*}
\hat U_n\left(x\right)&=\frac{\sqrt{n}}{\sqrt{{\rm I}(\hat\vartheta_n
  )}}\int_{-\infty }^{x }\dot \ell (\hat\vartheta_n ,y)\,\left[ {\rm
  d}  \hat F_n\left(y\right)-f(\hat\vartheta_n ,y){\rm d}y\right]\\
&=\frac{1}{\sqrt{{\rm I}(\hat\vartheta_n
  )}}\int_{-\infty }^{x }\dot \ell (\hat\vartheta_n ,y)\, {\rm
  d}  \sqrt{n}\left[\hat F_n\left(y\right)-F(\vartheta_0 ,y)\right]\\
&\quad +\frac{\sqrt{n}}{\sqrt{{\rm I}(\hat\vartheta_n
  )}}\int_{-\infty }^{x }\dot \ell (\hat\vartheta_n ,y)\,\left[
  f\left(\vartheta _0,y\right)-f(\hat\vartheta_n ,y)\right]{\rm d}y \\
&=\frac{1}{\sqrt{{\rm I}(\vartheta _0
  )}}\int_{-\infty }^{x }\dot \ell (\vartheta _0 ,y)\, {\rm
  d}  B_n\left(y\right)\\
&\quad -\frac{\sqrt{n}(\hat\vartheta_n-\vartheta _0)}{\sqrt{{\rm I}(\vartheta _0
  )}}\int_{-\infty }^{x }\dot \ell (\vartheta _0 ,y)\,\dot
  f(\vartheta _0 ,y){\rm d}y+o\left(1\right)\\
&\Longrightarrow \frac{1}{\sqrt{{\rm I}(\vartheta _0
  )}}\int_{-\infty }^{x }\dot \ell (\vartheta _0 ,y)\, {\rm
  d}  B\left(F\left(\vartheta _0,y\right) \right)\\
&\quad -\frac{1}{\sqrt{{\rm I}(\vartheta _0
  )}}\int_{-\infty }^{\infty  }\dot \ell (\vartheta _0 ,y)\, {\rm
  d}  B\left(F\left(\vartheta _0,y\right) \right)\;\int_{-\infty }^{x
  }\frac{\dot \ell (\vartheta _0 ,y)\,\dot
  f(\vartheta _0 ,y)}{{\rm I}(\vartheta _0
  )}{\rm d}y.
\end{align*}
Let us put $F\left(\vartheta_0 ,x\right)=t$, $F\left(\vartheta_0 ,y\right)=s$
and $h\left(\vartheta _0,s\right)=\dot \ell (\vartheta _0 ,y\left(s\right))$,
where $y\left(\vartheta_0 ,s\right)$ is solution $y$ of this equation
$F\left(\vartheta_0 ,y\right)=s $. Then the limit process can be written as
follows
\begin{align*}
U\left(t\right)&=\int_{0}^{t}h\left(\vartheta _0,s\right){\rm
  d}B\left(s\right)-\int_{0}^{1}h\left(\vartheta _0,s\right){\rm
  d}B\left(s\right)\;\int_{0}^{t}h\left(\vartheta _0,s\right)^2{\rm d}s
\\ &=\int_{0}^{t}h\left(\vartheta _0,s\right){\rm
  d}w\left(s\right)-\int_{0}^{1}h\left(\vartheta _0,s\right){\rm
  d}w\left(s\right)\;\int_{0}^{t}h\left(\vartheta _0,s\right)^2{\rm
  d}s\\ &\qquad -w\left(1\right)\int_{0}^{t}h\left(\vartheta _0,s\right){\rm
  d}s=W\left(\tau \right)-W\left(1\right)\tau
-w\left(1\right)\int_{0}^{t}h\left(\vartheta _0,s\right){\rm d}s.
\end{align*}
Therefore the limit
statistic is not free of distribution and this approach does not allow to
construct the ADF GoF test.


\begin{thebibliography}{1}
\bibitem{DK} Dachian, S. and Kutoyants, Yu.A. (2007)  On the goodness-of-fit
tests for some continuous time processes, in {\it Stat. Models
for Bio.- Tech. Systems}, F.Vonta {\it et al.} (Eds),
Boston,  395-413.
\bibitem{Dar} Darling, D. A. (1955) The Cram\'er-Smirnov test in the parametric
  case. {\it  Ann. Math. Statist}., 26, 1-20.


\bibitem{Durett} Durett, R. (1996)  {\sl Stochastic Calculus. A practical
introduction}. Boca Raton: CRC Press.

\bibitem{Fou92} Fournie, E.  (1992)  Un test de type Kolmogorov-Smirnov
 pour processus de diffusions ergodic. {\sl    Rapport de Recherche,
 {\bf 1696}},  INRIA, Sophia-Antipolis.

\bibitem{FV} Freidlin, M. I. and  Wentzell, A. D.(1998)  {\it Random
  Perturbations of
  Dynamical Systems.} 2nd Ed., Springer, N.Y.

\bibitem{HB93} Hwang, S.Y. and Basawa, I.V. (1993). Asymptotic optimal inference for a
class of nonlinear    time    series    models.    {\it Stochastic    Process
Appl.}    46,    91-113


\bibitem{IK01} Iacus S., Kutoyants Yu. A. (2001)   Semiparametric hypotheses
  testing for dynamical systems with small noise.  {\it Math. Methods
     Statist.} 10, 1, 105-120.

\bibitem{Kh81} Khmaladze, E. (1981) Martingale approach in the theory of
  goodness-of-fit tests. {\it Theory Probab. Appl. }, 26, 240-257.

\bibitem{KK} Kleptsyna, M.,    Kutoyants Yu. A. (2013) On asymptotically
  distribution  free  tests with parametric hypothesis
  for    ergodic diffusion processes. To appear in {\it
    Statist. Inference Stoch. Processes,}  (arXiv:1305.3382).


\bibitem{Kut94} Kutoyants, Yu.A. (1994) {\it Identification of Dynamical
  Systems with Small Noise,} Kluwer, Dordrecht.
\bibitem{Kut98} Kutoyants, Yu.A. (1998) {\it Statistical Inference for Spatial
Poisson Processes,} Springer, N.Y.
\bibitem{Kut04} Kutoyants, Yu.A. (2004) {\it Statistical Inference for Ergodic
Diffusion Processes,} Springer, London.
\bibitem{K11}   Kutoyants, Yu. A., (2011) {On goodness-of-fit tests for perturbed
dynamical systems.}  {\it J. Statist. Plann. Inference}, 141,  1655-1666.
\bibitem{Kut12b} Kutoyants, Yu. A., (2013) On asymptotic distribution of
  parame-\\ter free  tests for   ergodic diffusion processes. To appear in {\it
    Statist. Inference Stoch. Processes,} (arXiv:1302.1026).

\bibitem{Kut13a}   Kutoyants, Yu. A., (2013) {On ADF GoF tests for perturbed
dynamical systems.} Sumitted.
\bibitem{Kut13b}   Kutoyants, Yu. A., (2014) On ADF    goodness-of-fit tests
  for  stochastic processes. To appear in   {\it New Perspectives
    on Stochastic Modeling and Data Analysis, }
J. Bozeman, V. Girardin and C. H. Skiadas (Eds).

\bibitem {KZ} Kutoyants Y.A. and Zhou, L.(2013) On approximation of the
  backward stochastic differential equation. To appear in {\it
    J. Statist. Plann. Inference}, (arXiv:1305.3728).

\bibitem{LS} Liptser, R. and Shiryayev, A.N. {\it Statistics of Random
  Processes.} v. 2, 2-nd ed. Springer, N.Y., 2005.

\bibitem{MTP} Maglapheridze, N., Tsigroshvili,  Z. P. and  van Pul, M. (1998)
 Goodness-of-fit tests tests for parametric hypotheses on the distribution of
 point   processes, {\it Math. Methods. Statist.} 7, 60-77.

\bibitem{NN} Negri, I. and Nishiyama, Y. (2009) { Goodness of fit test for
 ergodic diffusion processes}.  {\it Ann. Inst. Statist. Math.}, 61, 919-928.

\bibitem{NZ}   Negri, I., and
Zhou, L. (2012) On goodness-of-fit testing for ergodic diffusion process with
  shift parameter. To appear in {\it Statist. Inference
    Stoch. Processes}, (arXiv:1203.6547).
\bibitem{NY96} Yoshida, N. (1996) Asymptotic expansions for perturbed systems on
  Wiener space: maximum likelihood estimators,  {\it J. Multivariate Analysis},
  57, 1-36.
\end{thebibliography}
\end{document}